\documentclass[mnsc,nonblindrev]{informs3} 

\OneAndAHalfSpacedXI 

\usepackage{endnotes}
\let\footnote=\endnote
\let\enotesize=\normalsize
\def\notesname{Endnotes}%
\def\makeenmark{$^{\theenmark}$}
\def\enoteformat{\rightskip0pt\leftskip0pt\parindent=1.75em
  \leavevmode\llap{\theenmark.\enskip}}
\usepackage{enumerate}
\usepackage{graphicx}

\usepackage{comment}
\usepackage[small, margin=1cm]{caption}

\usepackage{appendix}

\usepackage{svg}

\usepackage{pgfplots, tikz}
\usepackage{pgfplotstable}
\usetikzlibrary{patterns}
\usepackage{psfrag,algorithm,ifthen,multirow,mathrsfs}
\usepackage{color}
\definecolor{strcolor}{rgb}{0.6, 0.2, 0.6}
\definecolor{commentcolor}{rgb}{0.3125, 0.5, 0.3125}
\definecolor{keycol}{rgb}{0, 0, 1}

\usepackage{amssymb}
\usepackage{amsmath}
\usepackage{bbm}
\usepackage{lscape}
\usepackage{nicefrac}
\usepackage{enumitem}
\usepackage{calc}
\usepackage{dsfont}

\usepackage{url}
\usepackage{color}   
\usepackage{hyperref}
\hypersetup{
    colorlinks=true, 
    linktoc=all,     
    linkcolor=blue,  
}

\usepackage{bm}
\usepackage{algpseudocode}
\usepackage{scalerel}
\usepackage{soul}
\usepackage[figuresleft]{rotating}
\usepackage{subfig}
\usepackage{booktabs}
\usepackage{marvosym}

\usepackage{tikz}
\usetikzlibrary{shapes,arrows}
\tikzstyle{tank} = [rectangle, draw, fill=green!10,
text width=6em, text badly centered, node distance=2.5cm, inner sep=0pt]
\tikzstyle{waterbody} = [rectangle, draw, fill=blue!10,
text width=5em, text centered, minimum height=3em]
\tikzstyle{line} = [draw, very thick, color=black!20, -latex']
\tikzstyle{bline} = [draw, very thick, color=blue!40, -latex']
\tikzstyle{rline} = [draw, very thick, color=red!40, -latex']
\tikzstyle{plant} = [draw, ellipse,fill=red!10, node distance=2.5cm, text width = 5em, text centered, minimum height=3em]
\tikzstyle{dotted} = [rectangle, draw, dashed, fill=white!100, text width = 5em, text centered, minimum height=3em]

\tikzstyle{white1} = [rectangle,fill=white!20, text width=6em, text centered, node distance=2.5cm, minimum height=2em]
\tikzstyle{white2} = [rectangle,fill=white!20, node distance=2.5cm, minimum height=2em]
\tikzstyle{blackdot} = [circle, fill = black!90, text = white]

\usepackage[export]{adjustbox}

\newcommand{\mb}[1]{\mbox{\boldmath \ensuremath{#1}}}

\newcommand{\mbh}[1]{\mb{\hat{#1}}}

\newcommand{\eg}{\textit{e.g.}}
\newcommand{\ie}{\textit{i.e.}}

\newcommand{\gged}[1]{{\color{brown} #1}}
\newcommand{\tz}[1]{{\color{blue} #1}}

\def\blot{\quad \mbox{$\vcenter{ \vbox{ \hrule height.4pt
      \hbox{\vrule width.4pt height.9ex \kern.9ex \vrule width.4pt}
           \hrule height.4pt}}$}}

\usepackage{natbib}
 \bibpunct[, ]{(}{)}{,}{a}{}{,}%
 \def\bibfont{\small}%
\TheoremsNumberedThrough     
\ECRepeatTheorems

\EquationsNumberedThrough    

\usepackage{accents}

\begin{document}


 \RUNAUTHOR{Loke, Zhu and Zuo}

\RUNTITLE{Optimize-via-Predict}

\TITLE{Optimize-via-Predict \\ \large{Realizing out-of-sample optimality in data-driven optimization} }

\ARTICLEAUTHORS{%
    \AUTHOR{Gar Goei Loke}
	\AFF{Rotterdam School of Management, Erasmus University, Burgemeester Oudlaan 50, 3062PA Rotterdam, The Netherlands, \EMAIL{loke@rsm.nl}}
    
    \AUTHOR{Taozeng Zhu}
	\AFF{Institute of Supply Chain Analytics, Dongbei University of Finance and Economics, 116025 Dalian, China, \EMAIL{taozeng.zhu@gmail.com}}	
	
	\AUTHOR{Ruiting Zuo}
	\AFF{Financial Technology Thrust, Society Hub, The Hong Kong University of Science and Technology (Guangzhou), 510000 Guangzhou, China, \EMAIL{ruitingzuo@ust.hk}}
}

\ABSTRACT{We examine a stochastic formulation for data-driven optimization wherein the decision-maker is not privy to the true distribution, but has knowledge that it lies in some hypothesis set and possesses a historical data set, from which information about it can be gleaned. We define a prescriptive solution as a decision rule mapping such a data set to decisions. As there does not exist prescriptive solutions that are generalizable over the entire hypothesis set, we define out-of-sample optimality as a local average over a neighbourhood of hypotheses, and averaged over the sampling distribution. We prove sufficient conditions for local out-of-sample optimality, which reduces to functions of the sufficient statistic of the hypothesis family. We present an optimization problem that would solve for such an out-of-sample optimal solution, and does so efficiently by a combination of sampling and bisection search algorithms. Finally, we illustrate our model on the newsvendor model, and find strong performance when compared against alternatives in the literature. There are potential implications of our research on end-to-end learning and Bayesian optimization.}%

\KEYWORDS{Data-driven optimization, Prescriptive analytics, Sufficient statistics, Robust optimization, Stochastic optimization, Finite-sample optimality}

\maketitle

\section{Introduction}
Data-driven optimization has become increasingly relevant \citep{Hertog_Postek_2016_bridging}. 
It is usually represented as the stochastic program:
\begin{equation*}
    \min\limits_q \; \mathbb{E}_{d\sim D}[C(q, d)],
\end{equation*}
where the goal is to determine optimal decision $q$ to a cost function $C$ under the true distribution $D$ of the uncertainty $d$, that is unknown. The decision-maker possesses a data set of historical observations, from which information about $D$ can be gleaned.

In stochastic optimization, one approximates the expectation $\mathbb{E}_{d\sim D}$ using the data, what is termed as sample average approximations \citep{shapiro2021lectures}. Similarly, there are parametric \citep[such as the predict-then-optimize framework,][]{fisher2014demand} and non-parametric ways (such as kernel \citep{scott2015multivariate} and polynomial interpolations \citep{turner2021efficient}) to estimate the underlying distribution and the expectation \citep{Deng_Sen_2018learning}. While these methods are consistent (asymptotically convergent), errors can occur in smaller data samples, as the data set is not representative of true distribution $D$. In end-to-end learning, this small data regime is particular critical \citep{gupta2021small}. Moreover, errors from estimation transfer to and are amplified in the optimization, leading to the optimizer's curse \citep{smith2006optimizer}. While some have approached this from the angle of regret minimization \citep{Chen_Xie_2021regret,poursoltani2023risk}, the more dominant stream in the last two decades is robust optimization.  

\subsection*{Data-driven robust optimization}

Robust optimization assumes that the uncertainty lies within an \emph{uncertainty} set of possible manifestations, and seeks to be robust to it \citep[see, \emph{e.g.},][]{ben1998robust} by sacrificing a small degree of performance \citep{ben2000robust}, leading to good out-of-sample performance \citep{gotoh2021calibration}. In data-driven robust optimization, the uncertainty is the distribution itself, a form of distributionally robust optimization \citep{delage2010distributionally}. 
Earlier works use summary statistics to define the uncertainty \citep[such as moments,][]{wiesemann2014distributionally}. The contemporaneous approach declares a divergence measure and constructs the uncertainty around a neighbourhood of the empirical distribution \citep{Natarajan_Pachamanova_Sim2009constructing, Long2022rsatisficing}. Popular divergence measures include Kullback-Leibler divergence \citep{Hu_Hong_2013kullback} and Wasserstein distance \citep{Esfahani_Kuhn2018data, Blanchet_Murthy_2019quantifying}. \cite{bertsimas2018robust} also studies the robustification of sample average approximations in stochastic programming. 
More recently, works examine data-driven optimization from the statistical inference standpoint \citep{Duchi_Glynn_Namkoong_2016statistics}. \cite{bertsimas2018data} considers distributions not significantly different from the data under statistical tests. Defining the uncertainty around unknown parameters of a distributional family is seen more broadly in end-to-end learning \citep{Zhu_Xie_Sim_2022joint}. Robust optimization has also been established to learn under regularization \citep{Xu2010robustregre, bertsimas2018characterization}. Specific works also construct regularizations from data-driven uncertainty sets \citep{Gao_Chen_Kleywegt_2017_wasserstein}. 

In most cases, while bounds on the probability that the uncertainty set would not contain the true distribution can be proven \citep{sutter2020general}, in general, it is difficult to assert such claims under out-of-sample assumptions, except in special cases \citep{Ben-TalNemirovski00,bertsimas2004price}. We posit that the reason is because the data sample itself is an uncertainty under the sampling distribution, rendering the uncertainty set a random set. 

\subsection{Data set as the uncertainty}
Most works neglect that the data set arises from the sampling distribution. In stochastic optimization, one cannot generalize out-of-sample from the empirical distribution. This justified distributionally robust optimization, but the latter has not solved this problem either. Uncertainty sets are functions of the data set, thus random sets! In both cases, it is wishful that the models built under specific training sets are expected to apply universally to all possible testing data. 

\cite{Liyanage05}'s seminal work on the newsvendor problem both illustrates the shortcomings of ignoring the sampling distribution and proposes an interesting solution (see Example \ref{ex.opstat}. Further attempts to generalize this approach to general data-driven optimization \citep{Chu_Shanthikumar_Shen_2008solving,jia2022sufficient} have been met with limited progress and the broader community has, thus far, not picked up on the deep conceptual ideas behind the work. 

\subsection*{Contributions} 

We frame prescriptive analytics as a \emph{decision rule} mapping data sets (the uncertainty under the sampling distribution) to decisions, without \emph{a priori} knowledge of the true distribution. As there exist no function that is out-of-sample optimal with respect to every distribution in a given hypothesis set, the decision-maker seeks locally optimal solutions within a subset, termed a \emph{localization}. 
\begin{enumerate}[label = \alph*)]
    \item We prove sufficient conditions for out-of-sample optimality in data-driven optimization for a given localization -- they are functions of sufficient statistics (Theorem \ref{thm.suffissuff}).
    \item We write out an optimization problem that yields such an optimal solution (Theorems \ref{thm.restriction} and \ref{thm.ovp_optimality}). Under some conditions, it reduces to a bisection search, and is efficiently solved.
    \item We test our model on the newsvendor problem. It is superior to alternatives. We illustrate specificity-sensitivity trade-off in the selection of the localization.
\end{enumerate}

As our solution is a function of the maximum likelihood estimator, it bridges the idea that prescriptive analytics follows from predictive. Thus we term our approach \emph{optimize-via-predict}. Our work most closely relates to \cite{sutter2020general}; however, they consider asymptoptic optimality, whereas we focus on finite-sample optimality. Our work generalizes ideas in \cite{Liyanage05} to general decision policies and general convex problems.


\subsection*{Implications on end-to-end learning and other domains}

Contextual stochastic optimization, or end-to-end learning, assumes side information or contextual information $\bm x$ that helps with the inference of the uncertainty $d$ \citep{Bertsimas_Kallus_2020predictive}. It is sometimes viewed as a form of data-driven optimization, with structural assumptions \citep[conditional distribution with respect to context $\bm x$,][]{esteban2022distributionally}, and decisions $q$ are a function of context $\bm x$. This perspective is advocated by \cite{van2021data}. Nonetheless, methods are not restricted to data-driven optimization \citep[such as][]{Ban_Rudin_2019big,Elmachtoub_Grigas_2020smart}. By characterizing out-of-sample optimality for data-driven optimization, our work opens the door to potentially examine out-of-sample optimal end-to-end learning. 

Separately, the interpretation of the localization as prior exogenous information sets up the possibility of developing a Bayesian framework for prescriptive analytics \citep{Chu_Shanthikumar_Shen_2008solving}.

\smallskip
\noindent\textbf{Organization:} In \S\ref{sec.mod}, we define prescriptive solutions, localizations and out-of-sample optimality, culminating in \S\ref{subsec.suff} -- conditions for out-of-sample optimality and our proposed model. In \S\ref{sec.illust}, we illustrate on the newsvendor problem. 
Proofs of Propositions \ref{prop.upbound}, \ref{prop.suffstat} and \ref{prop.beth_conv} are omitted as they follow immediately from definitions or well-known facts.

\noindent\textbf{Notation:} Denote $\mathcal{M}(\Theta, \mathcal{X})$ as the set of probability distributions on support $\Theta$ and outcomes $\mathcal{X}$.

\section{Out-of-sample Optimality in Data-driven Optimization}\label{sec.mod}

Consider a cost function $C(q, d)$ that depends on a decision variable $q \in \mathcal{Q}$ (`quantity') in feasible set $\mathcal{Q}$, and an uncertain variable $d\in\mathcal{D}$ (`demand') modelled by p.d.f.s within a hypothesis set $\mathcal{H} := \{f(d;\theta) : \theta \in \Theta \subseteq \mathbb{R}^o\}$, containing the true distribution. The function $f$ is assumed to be given. 

\begin{assumption}[Convex costs]\label{assump.convex}
The feasible set for the decisions $\mathcal{Q}$ is convex; and for all $d\in\mathcal{D}$, the cost function $C$ is convex in $q$.
\end{assumption}

Given any decision $q\in\mathcal{Q}$ and true parameter $\theta$, the decision-maker incurs an expected cost of 
\begin{equation}\label{eq.pure_opt}
    \phi(q;\theta) := \int_{\mathcal{D}} C(q, d) f(d; \theta) \,dd.
\end{equation}

In reality, the decision-maker seeks to minimize expected costs $\phi(q; \theta)$ by optimizing $q$. However, without perfect information, \emph{i.e.}, knowledge of the unobserved parameter $\theta$, this function is uncertain. Instead, the decision-maker is armed with data set $\{y_n\in\mathcal{D}\}_{n=1}^N$, which we represent as a vector $\bm y$, containing $N$ i.i.d. data points sampled from the distribution $f(\cdot; \theta)$.  


\begin{definition}[Prescriptive solution]
    A \emph{prescriptive solution} is a function $q:\mathcal{D}^N \rightarrow \mathbb{R}$ that maps a data set of size $N$ to a quantity, $\bm y \mapsto q(\bm y)$, with no explicit dependence on parameter $\theta$.
\end{definition}

Though in practice, the decision-maker is only availed one data set $\bm y$, and thus only wishes to solve for a single quantity $q(\bm y)$, the data $\bm y$ is a random variable under the sampling distribution. Thus, to verify the effectiveness of $q$, one needs to average over the sampling distribution. Otherwise, the decision $q(\bm y)$ is non-generalizable. 

\begin{definition}\label{def.presc_sol}
Denote the sampling distribution of data sets of size $N$ as $\bm y \in \mathcal{D}^N:= \mathcal{Y}$, with joint distribution $g(\bm y; \theta) := \prod\limits_{n=1}^N f(y_n; \theta)$.
\begin{enumerate}[label= \roman*.]
    \item The \emph{out-of-sample performance} of a prescriptive solution $q(\cdot)$ is 
\begin{equation}\label{eq.data-driven}
    \Phi[q(\cdot) | \theta]:= \mathbb{E}_{\bm y}\left[ \phi(q(\bm y)) \middle| \theta  \right] = \int_{\mathcal{Y}}\int_{\mathcal{D}} C(q(\bm y), d) f(d; \theta) g(\bm y; \theta) \,dd \,d\mb{y};
\end{equation} 
    \item \label{def.with_h} Let the \emph{localization} $u\in\mathcal{M}(\Theta, \mathbb{R})$ be a density. The \emph{expected out-of-sample performance} of a prescriptive solution $q(\cdot)$ with respect to the localization $u(\cdot)$ is
    \begin{equation}\label{eq.exp_perf}
        \Psi[q(\cdot); u]:= \mathbb{E}_{u}\big[ \Phi[q(\cdot) | \theta ] \big] = \int_{\Theta}\int_{\mathcal{Y}}\int_{\mathcal{D}} C(q(\bm y), d) f(d; \theta) g(\bm y; \theta) u(\theta) \,dd \,d\mb{y} \, d\theta.
    \end{equation}
    We say that the solution $q(\cdot)$ is \emph{out-of-sample locally optimal} about $u$ if it minimizes \eqref{eq.exp_perf}.
\end{enumerate}
\end{definition}

\smallskip

A prescriptive function $q(\cdot)$ that minimizes the functional $\Phi[\cdot | \theta]$ would be optimal. However, (i) $\Phi$ is ill-defined without knowledge of true $\theta$, and (ii) in general, there does not exist a prescriptive function $q(\cdot)$, independent of $\theta$, that would be optimal for \eqref{eq.data-driven} for all $\theta\in\Theta$ (see Example \ref{ex.futile_pto}). 
These reasons motivate restricting to a neighbourhood of $\Theta$, characterized by density $u$ in definition \ref{def.with_h}. Also, by integrating over $\theta$, dependence of $q$ on true $\theta$ is removed. The localization $u$ represents our region of interest for $\theta$. The narrower $u$, the better the prescriptive solutions performs on $u$. This is reminiscent of the specificity versus sensitivity trade-off, and we return to that in \S\ref{sec.illust}.

Taking expectations over the sampling distribution $\bm y$ results in an out-of-sample metric. As it would be a mouthful to repeat the term `out-of-sample', we hereon drop it. 

\subsection{Motivating our approach from examples}


\begin{example}[Predict-then-optimize]
In predict-then-optimize (PTO), the decision-maker estimates $\hat\theta$, such as the maximum likelihood estimator (MLE), $\hat{\theta}_{\scriptscriptstyle{\text{MLE}}} := \argmin_{\theta} g(\bm y; \theta)$, and then chooses the prescriptive solution:
\begin{equation}\label{eq.pto}
    q_{\scriptscriptstyle{\text{PTO}}}(\bm y) := \argmin\limits_{q}\int_{\mathcal{D}} C(q, d) f(d; \hat{\theta}(\bm y)) \,dd, \qquad \forall \bm y \in \mathcal{Y}_N.
\end{equation} 
\end{example}

The decision-maker treats $\hat\theta_{\scriptscriptstyle{\text{MLE}}}$ as true $\theta$.  In reality, it changes with and inherits error from the data. These errors can be amplified within $\phi(\cdot)$ \citep{smith2006optimizer}, leading to over-fitting. 

\begin{example}\label{ex.futile_pto}
An alternative that considers $\hat\theta$ varying over the sampling distribution of $\bm y$, by optimizing \eqref{eq.data-driven} assuming that $\theta$ was in fact $\hat{\theta}(\bm y)$, is:
\begin{equation}\label{eq.pto_more}
    q(\bm y) = \argmin\limits_{q(\cdot)} \int_{\mathcal{Y}_N} \int_{\mathcal{D}} C(q(\bm y), d) f(d; \hat{\theta}(\bm y))\, g(\bm y; \hat{\theta}(\bm y)) \,dd \,d\mb{y},
\end{equation}
Unfortunately, it has the same optimality conditions as \eqref{eq.pto} and thus is no different from $q_{\scriptscriptstyle{\text{PTO}}}$. In other words, just because one considers the out-of-sample objective, it does not necessarily (i) lead to a different solution, nor (ii) yield a prescriptive solution that is independent of true $\theta$.
\end{example}

\begin{example}[\cite{Liyanage05}]\label{ex.opstat}
Here, $C$ is the newsvendor problem and $\mathcal{H}$ the family of exponential distributions. They propose the prescriptive solution, termed \emph{operational statistics},  $q_{\scriptscriptstyle{\text{OS}}}(\bm y) = \alpha \hat{\theta}(\bm y)$ for a specific constant $\alpha \in \mathbb{R}$ independent of $\bm y$ and true $\theta$, where $\hat{\theta}(\bm y)$ is the MLE. This solution dominates PTO, \emph{i.e.}, $\Phi[q_{\scriptscriptstyle{\text{OS}}} | \theta] \leq \Phi[q_{\scriptscriptstyle{\text{PTO}}} | \theta]$ for all $\theta\in\Theta$. If $\mathcal{H}$ is the set of empirical distributions, with $\hat\theta$ being the order statistic, their solution is again linear in $\hat\theta$. 
\end{example}

There are two important points of note here. First, while $q_{\scriptscriptstyle{\text{OS}}}$ dominates $q_{\scriptscriptstyle{\text{PTO}}}$ over every $\theta\in\Theta$, the operational statistic is not optimal for \eqref{eq.data-driven}. The following Example illustrates this. 
A more complicated function of the MLE can have stronger performance than $q_{\scriptscriptstyle{\text{OS}}}$, at least over a subset of $\Theta$. 
This supports ideas in Definition \ref{def.presc_sol} -- by restricting to a localization, one could potentially obtain a prescriptive solution with a certificate of local out-of-sample optimality.

\begin{figure}[!h]
    \centering
    \includegraphics[width = 10cm]{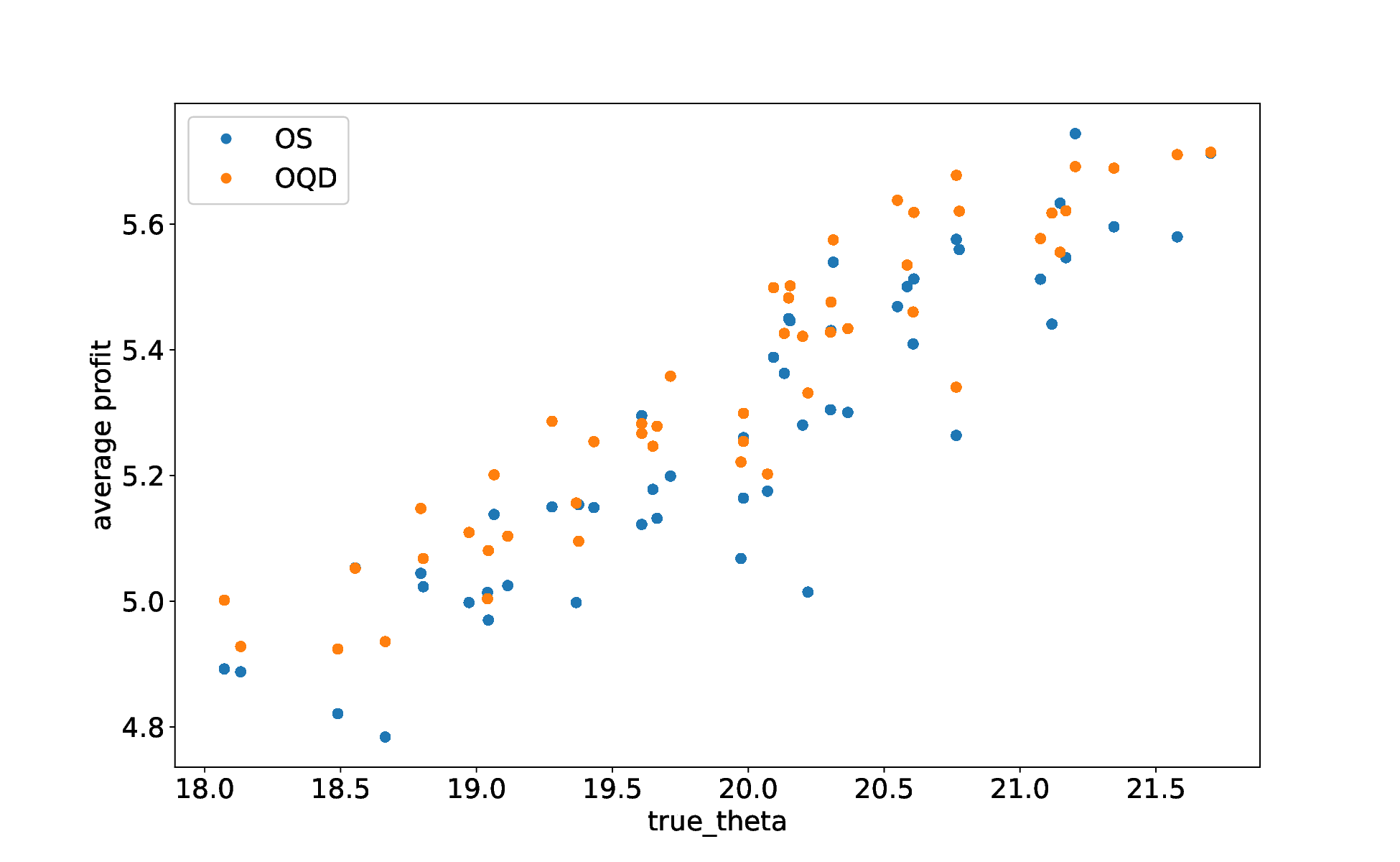}
    \caption{Performance of operational statistics (OS -- blue) against quadratic variant (OQD -- yellow)}
    \label{fig.quad_compare}
\end{figure}

\begin{example}\label{example.quad}
    Consider the prescriptive solution $q_{\scriptscriptstyle{\text{OQD}}}(\bm y) := \alpha\hat\theta(\bm y) - \hat\theta(\bm y)^2/2N^{3} $. Figure \ref{fig.quad_compare} shows that $q_{\scriptscriptstyle{\text{OQD}}}$ can dominate $q_{\scriptscriptstyle{\text{OS}}}$ in some region of true $\theta$. 
\end{example} 


\subsection{Sufficient statistics are sufficient}\label{subsec.suff}

The important observation from Example \ref{ex.opstat} is that the mean and order statistics are sufficient (minimal) statistics for the exponential and empirical distributions; and the decision is a linear decision rule of them. In Example \ref{example.quad}, we considered a quadratic decision rule. This observation was noted by \cite{jia2022sufficient}, without means of exploiting it. Here, we propose prescriptive solutions constructed from sufficient statistics. This extends ideas behind $q_{\scriptscriptstyle{\text{OS}}}$ to general candidate policies and convex optimization problems. Interestingly, sufficient statistics were employed in proofs of out-of-sample optimality \citep[such as][]{sutter2020general}, confirming suspicions of their role.

\begin{definition}[MLE-Sufficient Prescriptive Solution] {\color{white} a}
    \begin{enumerate}[label = \roman*.]
    \item A prescriptive solution $q(\bm y)$ is \emph{MLE-sufficient} if it can be written as $q(\hat\theta)$, with only explicit dependence on the MLE $\hat\theta(\bm y)$. Denote the class of MLE-sufficient prescriptive solutions as $\mathcal{Q}^{\mathcal{S}}$.

    \item For a given localization $u$, a candidate prescriptive solution $q_{\scriptscriptstyle{\text{OVP}}}$ to the optimization problem
    \begin{equation}\label{eq.ovp} \tag{OVP}
        \inf\limits_{q(\cdot)\in\mathcal{Q}^{\mathcal{S}}} \Psi[q(\cdot);u],
    \end{equation}
    is called the \emph{optimize-via-predict} (OVP) solution with respect to localization $u$.
    \end{enumerate}
\end{definition}

\begin{proposition}\label{prop.upbound}
    Any OVP solution forms an upper bound to minimizing \eqref{eq.exp_perf}.
\end{proposition}

\begin{proposition}\label{prop.suffstat}
    Let $\hat\theta$ be a sufficient statistic for $f(\cdot, \theta)$. There exists non-negative functions $g_0$ and $g_1$ such that the joint density decomposes: $g(\bm y; \theta) = g_0(\bm y) g_1(\hat{\theta}(\bm y), \theta), \;\;\forall \bm y\in\mathcal{Y}, \forall \theta\in\Theta$.
\end{proposition}

Our intent is to transfer the problem from $\bm y \in \mathcal{Y}$ onto $\hat\theta \in \Theta$, which is of far smaller dimension.

\begin{assumption}[Non-negative Jacobian]\label{assump.jacob}
There is a re-parameterization from the space $\bm y \in \mathcal{Y}$ to $(\hat\theta, \bm y|\hat\theta) \in \Theta \times \mathcal{D}^{N-o}$, $\Theta \subseteq \mathbb{R}^{o}$, where we abused the notation $\bm y|\hat\theta$ to refer to some parametrization of the restriction $\{\bm y: \hat{\theta}(\bm y) = \hat{\theta} \}$ for a given $\hat\theta$. We further assume that the Jacobian for the change of variables, denoted as $J(\hat{\theta}, \bm y | \hat{\theta})$, exists and is non-negative for all $\hat\theta \in \Theta$ and $\bm y | \hat\theta \in \mathcal{D}^{N-o}$. 
\end{assumption}

\begin{theorem}[Restricted problem]\label{thm.restriction}
If we restrict the search for prescriptive solutions to $\mathcal{Q^S}$, the set of MLE-sufficient prescriptive solutions, then \eqref{eq.ovp} reduces to
\begin{equation} \label{eq.reform}\tag{R-OVP}
    \inf\limits_{q(\cdot)\in\mathcal{Q^S}} \int_{\bm \Theta} K(\hat\theta) \, \beth(q, \hat\theta) \,d\hat{\theta}, 
\end{equation}
where \begin{equation}\label{eq.K}
    K(\hat{\theta}) := \int_{\{\bm y: \hat{\theta}(\bm y) = \hat{\theta} \}} g_0(\hat\theta, \bm y | \hat\theta) J(\hat{\theta}, \bm y | \hat{\theta}) \, d(\bm y|\hat{\theta}) \quad \mbox{and} \quad \beth(q, \hat{\theta}) = \int_{\bm \Theta} \phi(q, \theta) g_1(\hat{\theta}, \theta) u(\theta) \, d\theta,
\end{equation}
and $g_0$ and $g_1$ are defined in Proposition \ref{prop.suffstat}.
\end{theorem}

\begin{proof}{Proof of Theorem \ref{thm.restriction}.} Follows from Proposition \ref{prop.suffstat} and the change of variables $\bm y \rightarrow (\hat\theta, \bm y|\hat\theta)$. \hfill\Halmos
\end{proof}



\begin{theorem}[Sufficiency]\label{thm.suffissuff}
    Under Assumptions \ref{assump.convex} and \ref{assump.jacob}, there exists some MLE-sufficient prescriptive solution $\tilde{q}(\cdot) \in \mathcal{Q}^{\mathcal{S}}$ that is out-of-sample locally optimal for a given localization $u$. 
\end{theorem}

\begin{proof}{Proof of Theorem \ref{thm.suffissuff}.}
Let $K(\hat\theta)$ be defined in the sense of \eqref{eq.K}. Define for every given $\hat\theta\in\Theta$,
\begin{equation}\label{eq.weird_dens}
    \mathscr{X}(\hat\theta, \bm y^\circ|\hat\theta) := \frac{g_0(\hat\theta, \bm y^\circ | \hat\theta) J(\hat{\theta}, \bm y^\circ | \hat{\theta})}{K(\hat\theta)},
\end{equation}
so that $\mathscr{X}$ is a density, as by Proposition \ref{prop.suffstat} and Assumption \ref{assump.jacob}, $ \mathscr{X}(\hat\theta, \bm y^\circ|\hat\theta) \geq 0 $ and,
\begin{equation*}
     \int_{\{\bm y^\circ: \hat{\theta}(\bm y^\circ) = \hat{\theta} \}} \mathscr{X}(\hat\theta, \bm y^\circ|\hat\theta) \, d(\bm y^\circ|\hat{\theta}) = 1.
\end{equation*}
Given an optimal prescriptive solution $q(\cdot)$ for \eqref{eq.ovp}, construct a new solution $\tilde{q}(\cdot)$ as follows:
\begin{equation}\label{eq.construction}
    \tilde{q}(\bm y) := \int_{\{\bm y^\circ: \hat{\theta}(\bm y^\circ) = \hat{\theta} \}} q(\hat\theta, \bm y^{\circ}|\hat\theta) \mathscr{X}(\hat\theta, \bm y^\circ|\hat\theta) \, d(\bm y^\circ|\hat{\theta}).
\end{equation}
By construction, $\tilde{q}(\cdot)$ is only explicitly in $\hat\theta$, and thus is MLE-sufficient. It is feasible for all $\bm y$ due to convexity of the decision space $\mathcal{Q}$ and that $\mathscr{X}$ is a density. Its objective value is
\begin{align}
    & \int_{\Theta}\int_{\mathcal{Y}} \phi(\tilde{q}(\bm y), \theta) g(\bm y; \theta) u(\theta) \,d\mb{y} \, d\theta, \nonumber \\
    = & \int_{\Theta} \int_{\Theta} \int_{\{\bm y: \hat{\theta}(\bm y) = \hat{\theta} \}} g_0(\hat\theta, \bm y | \hat\theta) J(\hat\theta, \bm y | \hat\theta) g_1(\theta, \hat\theta) u(\theta) \cdot \nonumber \\
    & \quad  \phi\left(\int_{\{\bm y^\circ: \hat{\theta}(\bm y^\circ) = \hat{\theta} \}} q(\hat\theta, \bm y^{\circ}|\hat\theta) \mathscr{X}(\hat\theta, \bm y^\circ|\hat\theta) \, d(\bm y^\circ|\hat{\theta}) \; , \; \theta \right) d\bm y|\hat\theta \, d\hat\theta \, d\theta \nonumber \\
    \leq & \int_{\Theta} \int_{\Theta} \int_{\{\bm y: \hat{\theta}(\bm y) = \hat{\theta} \}} g_0(\hat\theta, \bm y | \hat\theta) J(\hat\theta, \bm y | \hat\theta) g_1(\theta, \hat\theta) u(\theta) \cdot \label{eq.thmconv} \\
    & \quad \left( \int_{\{\bm y^\circ: \hat{\theta}(\bm y^\circ) = \hat{\theta} \}} \phi \left(q(\hat\theta, \bm y^{\circ}|\hat\theta), \theta \right) \mathscr{X}(\hat\theta, \bm y^\circ|\hat\theta) \, d(\bm y^\circ|\hat{\theta}) \right) d\bm y|\hat\theta \, d\hat\theta \,d\theta \nonumber\\
    = & \int_{\Theta} \int_{\Theta} \left( \int_{\{\bm y: \hat{\theta}(\bm y) = \hat{\theta} \}} g_0(\hat\theta, \bm y | \hat\theta) J(\hat\theta, \bm y | \hat\theta) d\bm y|\hat\theta \right) g_1(\theta, \hat\theta) u(\theta) \cdot \label{eq.thmsup}\\
    & \quad \left( \int_{\{\bm y^\circ: \hat{\theta}(\bm y^\circ) = \hat{\theta} \}} \phi \left(q(\hat\theta, \bm y^{\circ}|\hat\theta), \theta \right) \mathscr{X}(\hat\theta, \bm y^\circ|\hat\theta) \, d(\bm y^\circ|\hat{\theta}) \right) d\hat\theta \,d\theta \nonumber \\
    = & \int_{\Theta} \int_{\Theta} K(\hat\theta) g_1(\theta, \hat\theta) u(\theta) \int_{\{\bm y^\circ: \hat{\theta}(\bm y^\circ) = \hat{\theta} \}} \phi \left(q(\hat\theta, \bm y^{\circ}|\hat\theta), \theta \right) \mathscr{X}(\hat\theta, \bm y^\circ|\hat\theta) \, d(\bm y^\circ|\hat{\theta}) d\hat\theta \,d\theta \nonumber \\
    = & \int_{\Theta} \int_{\Theta} \int_{\{\bm y^\circ: \hat{\theta}(\bm y^\circ) = \hat{\theta} \}} g_1(\theta, \hat\theta) u(\theta) \phi \left(q(\hat\theta, \bm y^{\circ}|\hat\theta), \theta \right) g_0(\hat\theta, \bm y^\circ|\hat\theta) J(\hat\theta, \bm y^\circ|\hat\theta) \, d(\bm y^\circ|\hat{\theta}) d\hat\theta \,d\theta \nonumber \\ 
    = & \int_{\Theta} \int_{\mathcal{Y}} \phi \left(q(\bm y), \theta\right) g(\bm y; \theta) u(\theta) \, d\bm y \, d\theta \label{eq.thmfinal}
\end{align}
where \eqref{eq.thmconv} follows from (i) Jensen's equality applied on the first argument of $\phi$, (ii) the construction of $\mathscr{X}$ as a density on $\bm y^{\circ}|\hat\theta$, and (iii) the non-negativity of $g_0$, $g_1$, $J$ and $u$; \eqref{eq.thmsup} holds as $g_0$ and $J$ are the only functions explicitly dependent on $\bm y|\hat\theta$; and till \eqref{eq.thmfinal}, we apply the definitions of $K$ and Proposition \ref{prop.suffstat}.
But \eqref{eq.thmfinal} is the objective value of prescriptive solution $q(\cdot)$. Hence, we have an MLE-sufficient prescriptive solution $\tilde{q}(\cdot)$ dominating the original one. Thus $\tilde{q}(\cdot)$ is optimal. \hfill\Halmos
\end{proof}

\begin{remark}
\begin{enumerate}[label = \alph*)]
    \item In the proofs of Theorems \ref{thm.restriction} and \ref{thm.suffissuff}, we did not use any property of the MLE, save for it being sufficient. Thus, they hold true for any sufficient statistic $\hat\theta$. The MLE being the minimal statistic, however, leads to the most succinct representation for $q$.
    \item If in the proof of Theorem \ref{thm.suffissuff}, $q$ is MLE-sufficient, then $q$ has no component in $\bm y | \hat\theta$ and thus it makes $\phi$ constant (hence, linear) in the argument of $\bm y|\hat\theta$. This meets the equality conditions for Jensen's inequality, thus the construction $\tilde{q}$ would not lead to a strictly better solution.
    \item If $\mathcal{H} = \mathcal{M}(\mathbb{R}, \mathbb{R})$, that is, that there are no distributional assumptions, then the order statistic is sufficient. However, it has dimensions $N$, so no reductions in dimensional complexity in $q$ is obtained, though structurally we obtain a decision rule that is a function of the order statistic.
\end{enumerate}
    
\end{remark}

Theorem \ref{thm.suffissuff} shows that we should seek MLE-sufficient prescriptive solutions. 
We explain the intuition behind these results. A sufficient statistic `fully captures all of the information about the distribution'. It thus seems natural that any good solution must also contain all of this information and thus be a function of some sufficient statistic. 
Indeed, the proof of Theorem \ref{thm.suffissuff} is a proof by symmetry. Directions away from the subspace spanned by the sufficient statistic are averaged away by the construction \eqref{eq.construction}, with some smart reweighing $\mathscr{X}$, leaving terms that are symmetric about the sufficient statistic. In other words, directions outside of the sufficient statistic are irrelevant. 

\subsection{Solving for locally optimal solutions}

\begin{theorem}[Optimality conditions for OVP]\label{thm.ovp_optimality}
    Suppose $q_{\scriptscriptstyle{\text{OVP}}}$ is a solution that point-wise minimizes $\beth(\cdot, \hat\theta)$ for every $\hat\theta \in \Theta$. Then it is locally optimal. 
\end{theorem}
\begin{proof}{Proof of Theorem \ref{thm.ovp_optimality}.}
    Suppose $q_{\scriptscriptstyle{\text{OVP}}}$ is not optimal for \eqref{eq.reform}. Hence, there exists some other precriptive solution $Q$ such that 
    \begin{equation}\label{eq.contradiction}
    \int_{\bm \Theta} K(\hat\theta) \, \beth(Q(\hat \theta) , \hat\theta) \,d\hat{\theta} - \int_{\bm \Theta} K(\hat \theta) \, \beth(q_{\scriptscriptstyle{\text{OVP}}}(\hat \theta), \hat\theta) \,d\hat{\theta} < 0.    
    \end{equation}
    However, by change of variables from $(\hat\theta, \bm y | \hat\theta)$ back to $\bm y$, we obtain that the LHS of the above is 
    \begin{equation*}
        \int_{\mathcal{Y}}  \left[\beth(Q(\hat \theta (\bm y)), \hat\theta (\bm y))  - \beth(q^*(\hat \theta (\bm y)), \hat\theta (\bm y)) \right]g_0(\bm y)  d\bm y \geq 0,
    \end{equation*}
    by non-negativity of $g_0$ and the optimality of $q_{\scriptscriptstyle{\text{OVP}}}$ for $\beth$ for all $\hat\theta\in\Theta$, contradicting \eqref{eq.contradiction}. \hfill\Halmos
    \end{proof}

Critically, $K$, which involves a high-dimensional integral, a non-standard domain set, an unwieldy parametrization $\bm y| \hat\theta$, and a Jacobian $J$, is not involved. Intuitively, as the sufficient statistic already `contains all necessary information', the information along $\bm y|\hat\theta$ can be safely discarded.

\begin{proposition}\label{prop.beth_conv}
    The function $\beth$ is convex in the first argument.
\end{proposition}

\begin{corollary}[Bisection search]\label{coro.bisection}
Under assumptions of Theorems \ref{thm.suffissuff} and \ref{thm.ovp_optimality},
    \begin{enumerate}[label = \roman*.]
    \item If $\beth(\cdot, \hat\theta)$ is differentiable for every $\hat\theta \in \Theta$, then $q_{\scriptscriptstyle{\text{OVP}}}$ is locally optimal if it satisfies
    \begin{equation}\label{eq.deriv=0}
    \frac{\partial \beth}{\partial q}(q,\hat\theta) = 0, \qquad \forall \hat\theta \in \Theta.
    \end{equation}
    \item If furthermore, $\phi(\cdot, \theta)$ is differentiable for every $\theta \in \Theta$, \eqref{eq.deriv=0} is equivalent to the condition,
    \begin{equation}\label{eq.ovp_optimality}
        \int_{\bm \Theta} \frac{\partial \phi}{\partial q} (q, \theta) g_1(\hat\theta, \theta) u(\theta) \, d\theta = 0, \qquad \forall \hat\theta \in \Theta.
    \end{equation}
    Moreover, the LHS of \eqref{eq.ovp_optimality} is monotone and thus, solving \eqref{eq.ovp_optimality} reduces to a bisection search problem on $q(\hat\theta)$ for every $\hat\theta \in \Theta$.
\end{enumerate} 
\end{corollary}

\begin{proof}{Proof of Corollary \ref{coro.bisection}.} (i) are first order conditions of Theorem \ref{thm.ovp_optimality}; in (ii), as $\phi$ is convex, $\partial\phi/\partial q$ is monotone. Also, $g_1$ and $u$ are non-negative. \hfill\Halmos
\end{proof}


In practice, we obtain one $\hat\theta$ for each data set, and thus when given a training data set, one only needs to evaluate $q(\cdot)$ at one point. Theorem \ref{thm.ovp_optimality} guarantees we can easily do that, as the optimality condition is a point-wise one. In the worst case, one only needs to perform a golden search, as Proposition \ref{prop.beth_conv} guarantees that $\beth$ is convex. We discuss computational strategies in Appendix \ref{append.codes}.

\medskip

\section{Illustration on the newsvendor problem}\label{sec.illust}

The newsvendor problem with selling and cost prices $p$ and $c$ respectively has profit, $R(q,d) =  p \min\{d, q\} - qc := -C(q,d)$, with order quantity $q$ and random demand $d$. If $f(\cdot, \theta)$ is exponentially-distributed with mean $\theta$, we can explicitly compute $\phi$ as follows: $\phi(q, \theta) = qc - p\theta\left(1- e^{-q/\theta}\right)$ and $\frac{\partial\phi}{\partial q}(q, \theta) = c - pe^{-q/\theta}$. $\phi$ is indeed differentiable, fulfilling Corollary \ref{coro.bisection}. We can approximate \eqref{eq.ovp_optimality} with a set of samples $\mathcal{U}^{M}:= \{\theta_m\}_{m=1}^{M}$ for the localization $u$, and solve for its (unique) zeroes
pointwise for any sample parameter estimate $\hat\theta(\bm y)$. A sample algorithm is provided in Algorithm \ref{algo.news}. We avoid diving into details about the simulation set up; the reader is referred to Appendix \ref{append.setup}.

\smallskip
\noindent\ul{Benchmarks}: We consider a total of seven benchmarks. First, we consider the two predict-then-optimize benchmarks, namely \eqref{eq.pto} and a sample average approximation (SAA) of \eqref{eq.pure_opt}. Second, we consider $q_{\scriptscriptstyle{\text{OS}}}$ in Example \ref{ex.opstat}, specifically the parametric version. Third, we consider four robust optimization models, namely, a vanilla robust optimization model with uncertainty on the unknown parameter $\theta$, and three DRO models with the moments, Wasserstein and KL-divergence uncertainty sets on unknown demand $d$. Specific formulations are available in Appendix \ref{append.benchmarks}. Calibration of radii are presented in Appendix \ref{append.cross}. We would have liked to examine DRO models with uncertain $\theta$. In practice, there is only one data set -- a single observation of $\theta$. Only in the case of the moments uncertainty set, it is possible to estimate the variance of $\theta$ using the sample variance of the sample mean. However, this leads to a nonconvex formulation, thus it is not considered as a benchmark. 

\subsection{Experiment 1: No misspecification}\label{subsec.no_misspec}
\begin{figure}[h] 
    \centering
    \vspace{-1.5cm}
    {{\includegraphics[height=10cm]{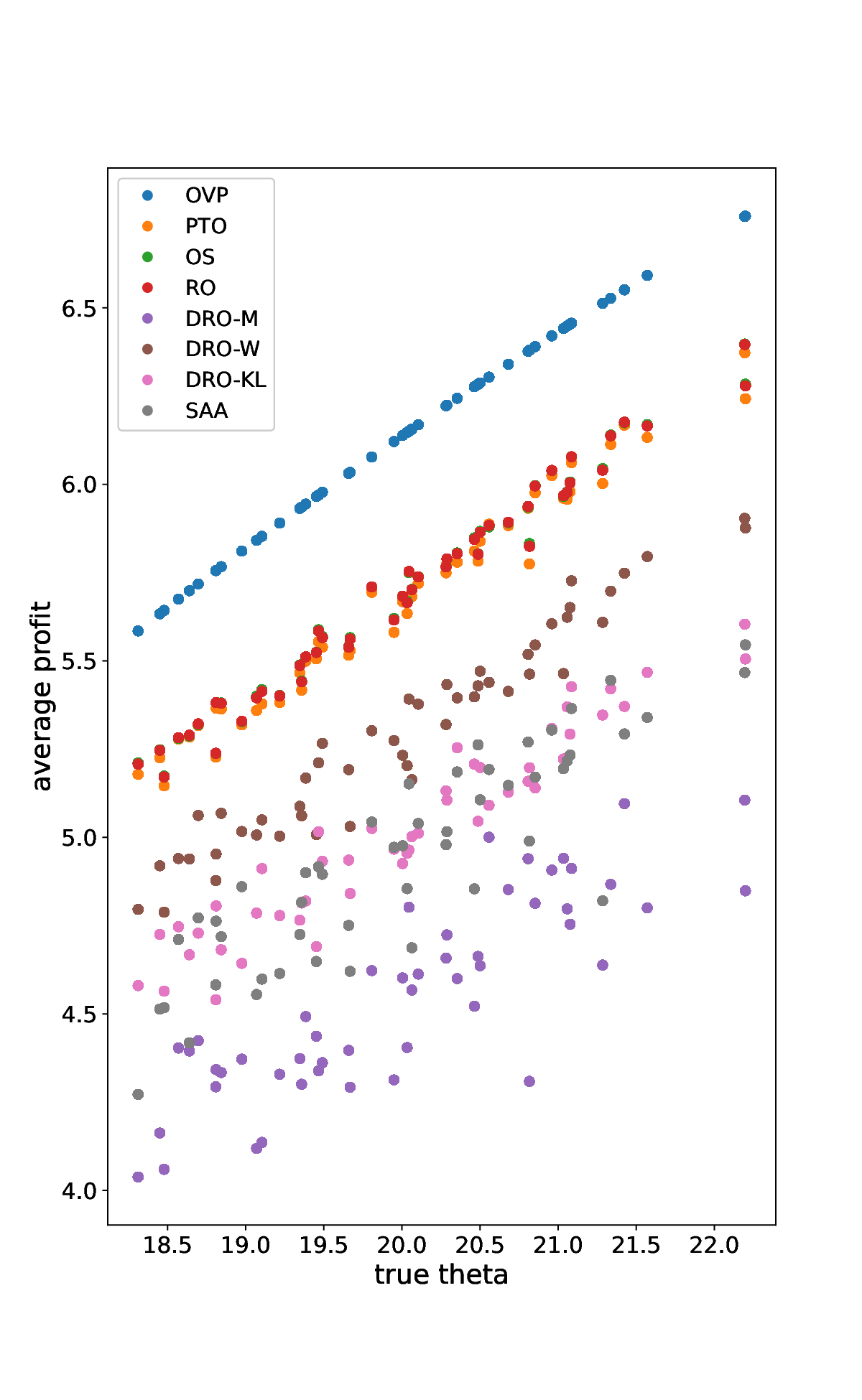}}}
   {{\includegraphics[height=10cm]{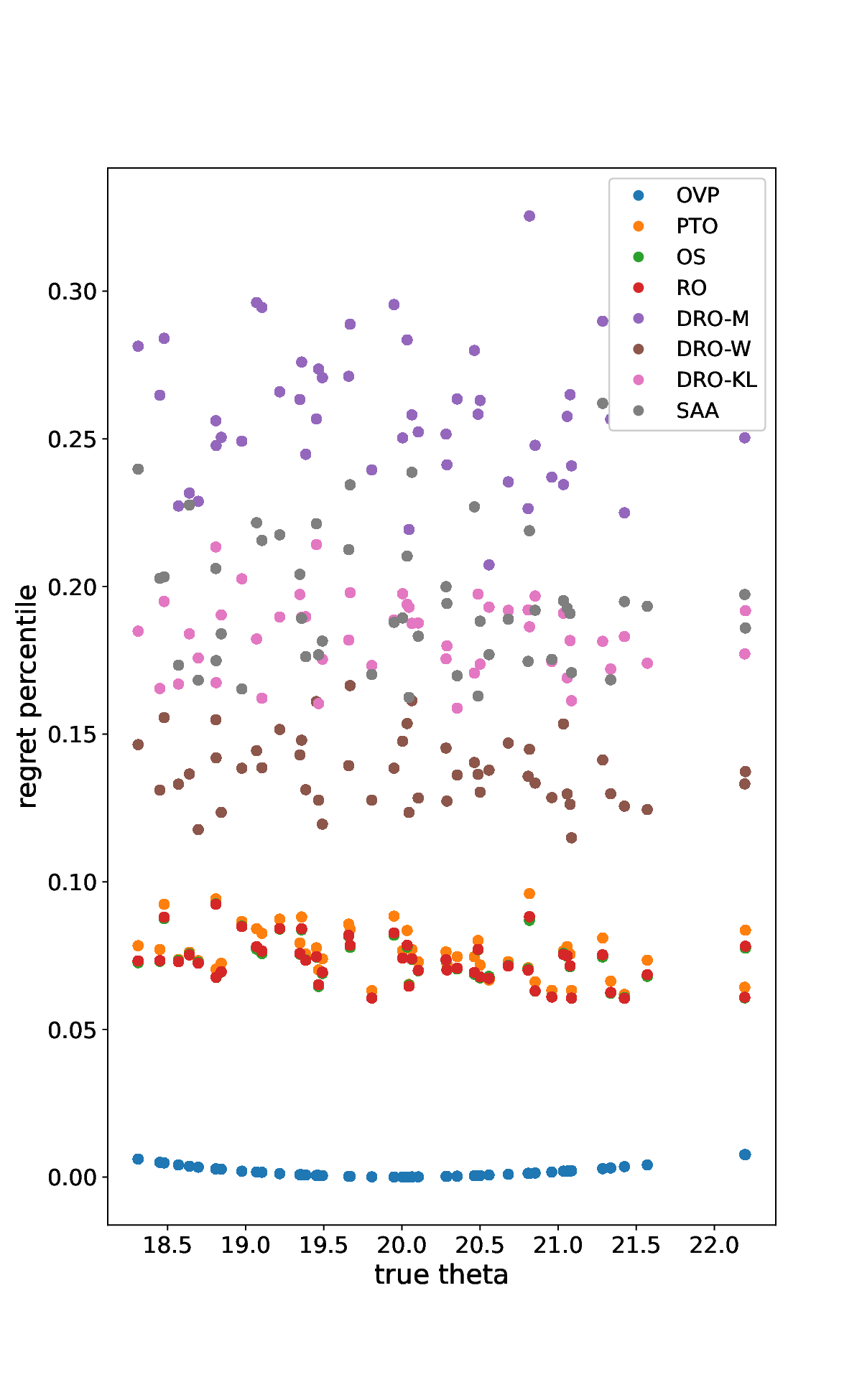} }}%
    \vspace{-0.8cm}
    \caption{Average profit (left) and ex-ante percentage regret (right) of different models over the range of $\theta$ in the localization $u\sim\mathcal{N}(20,1)$ }\label{fig.norm_base}%
    \vspace{-0.3cm}
\end{figure}

We first consider the case where the true distribution is indeed Exponential, \ie, there is no misspecification. In the base case, we consider the localization $u\sim\mathcal{N}(20,1)$. Figure \ref{fig.norm_base} plots the models' performance in terms of average profit, as well as percentage regret against the perfect information ex-ante oracle, which is just the quantile solution $F^{-1}(\frac{c}{p})$ assuming true $\theta$.

Most noticeably, OVP both outperforms the next best benchmark by a significant margin uniformly over the range of the localization, and achieves close to $0$ regret against the ex-ante oracle. What is interesting is the smoothness of the profit function for OVP, which results from OVP directly optimizing out-of-sample profit, as opposed to the benchmarks which use in-sample objective functions, thus affected by variations in the sampling distribution.

Amongst the benchmarks, models that explicitly assume an exponential distribution for the demand (PTO, OS [obscured by PTO and RO], RO and OVP) outperform those that do not (SAA and the DRO models). Among the latter, DRO-moments performs the worst, as its worst-case distributions are unlikely exponential \citep[if both mean and variance are exactly specified, it is a two-point distribution,][]{scarf1958min}. The data-driven RO models, being anchored on the data sample, result in worst-case distributions closer to an exponential distribution, thus outperforming the moments model. However, as their uncertainty set permits distributions far more diverse than the RO model, they pay the corresponding \emph{price of robustness}. The Wasserstein model outperforms SAA, being able to correct for sampling error. 
The KL divergence model obtained a different solution from SAA, but their out-of-sample performance are similar. 

Amongst exponential demand models, PTO is the worst. RO corrects for potential estimation error in $\theta$ and consequently outperforms PTO, though if its uncertainty set is too large, its performance will deteriorate. As proven by \cite{Liyanage05}, OS dominates PTO for all $\theta$, except the off-chance under extremes of the sampling distribution. Notably, what we gained using linear decision rules (\ie, the gap between PTO and OS) is only a small fraction of the gains from a general decision rule (\ie, the gap between OS and OVP). 

\begin{figure}[h] 
    \centering
    \vspace{-1.5cm}
   {{\includegraphics[height=10cm]{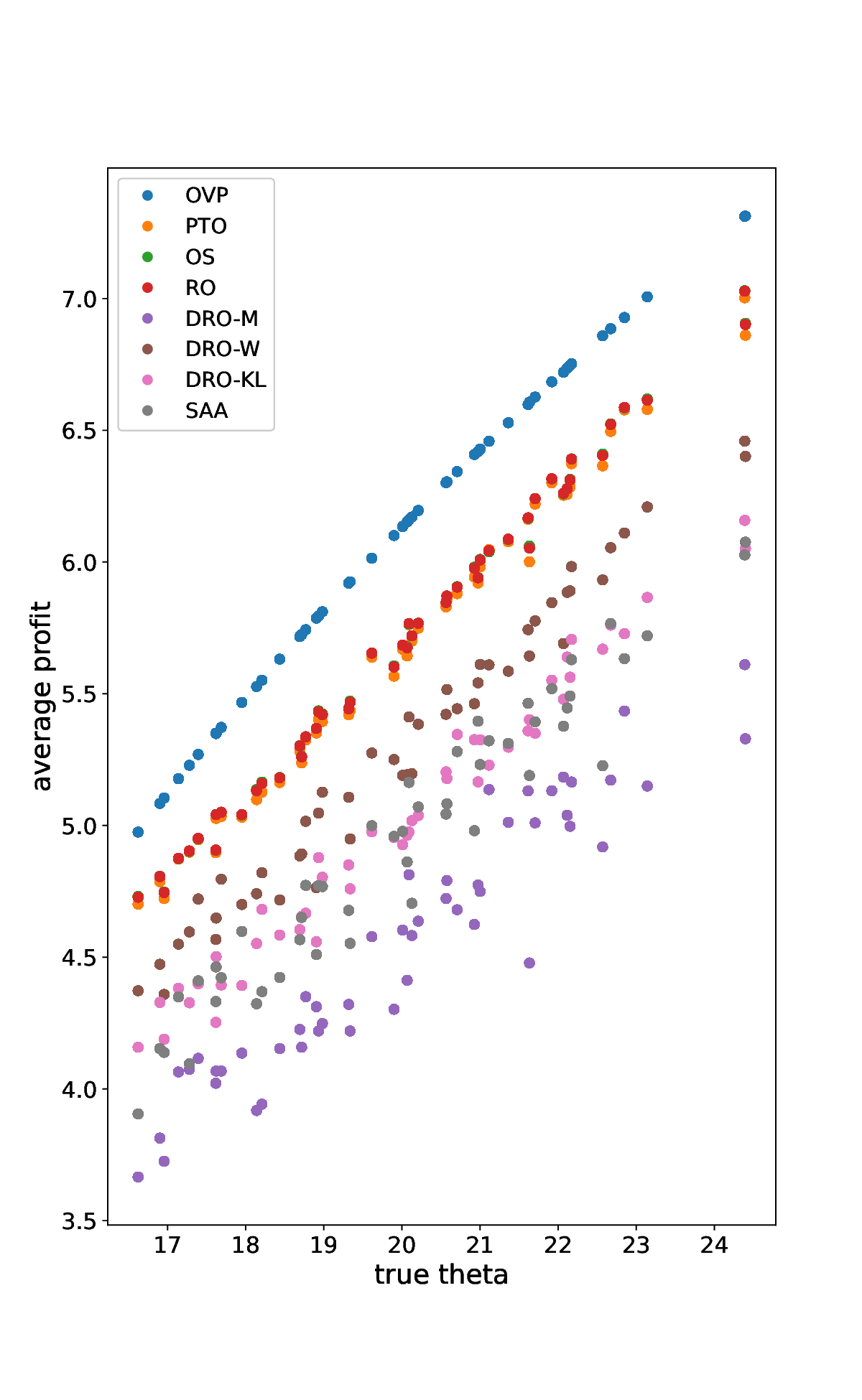}}}
   {{\includegraphics[height=10cm]{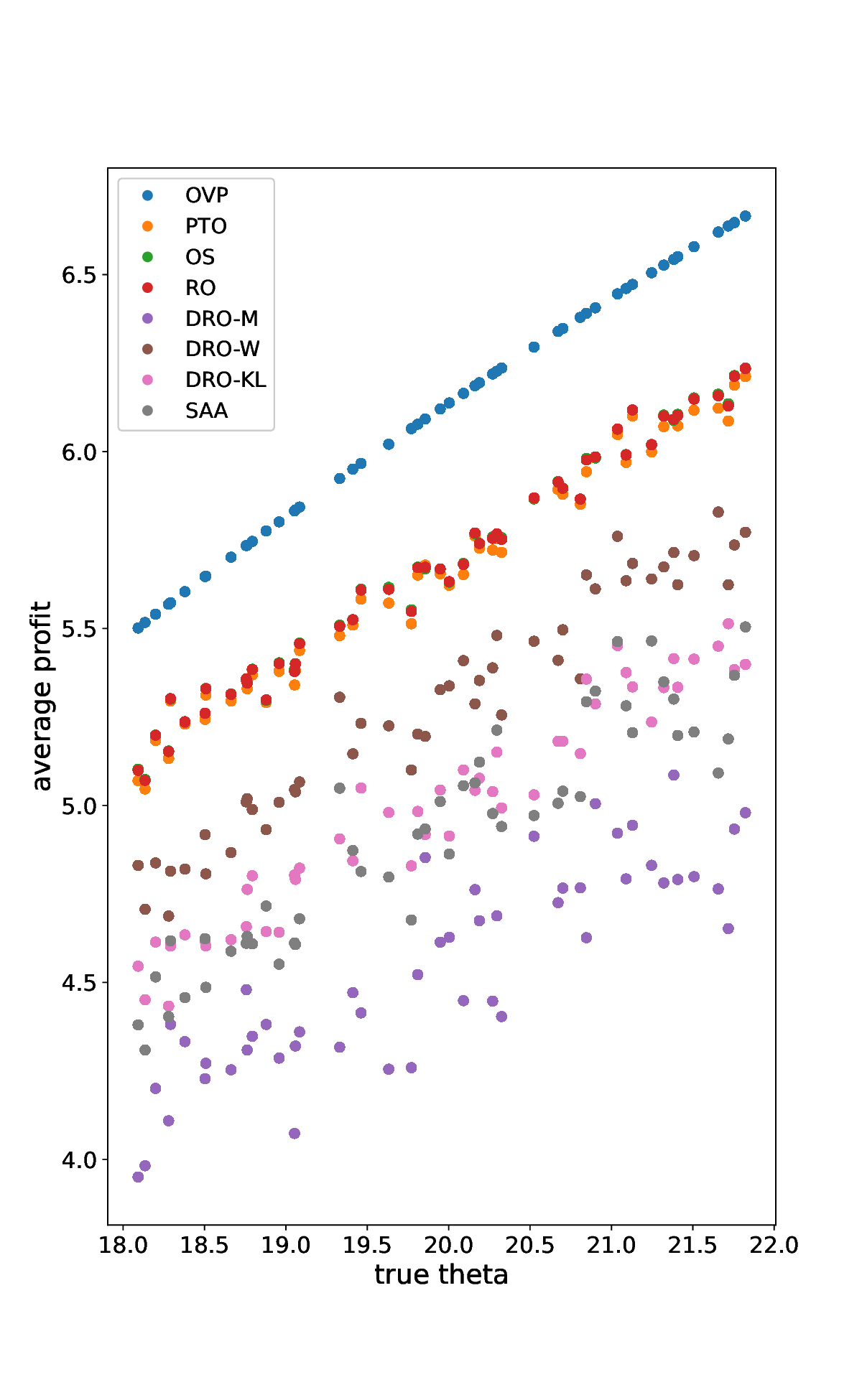} }}\\
    \vspace{-1.5cm}
    {{\includegraphics[height=10cm]{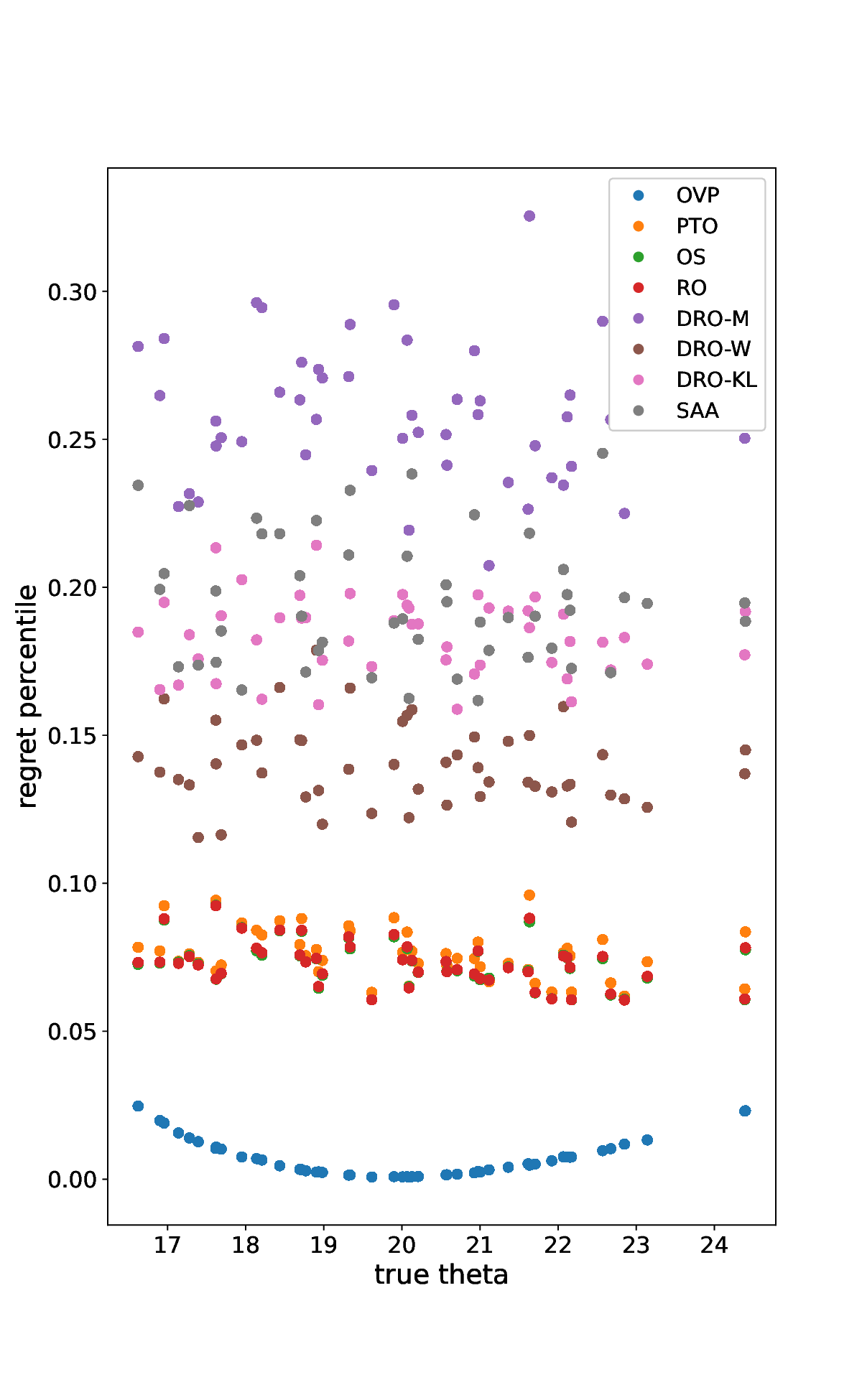}}}
    {{\includegraphics[height=10cm]{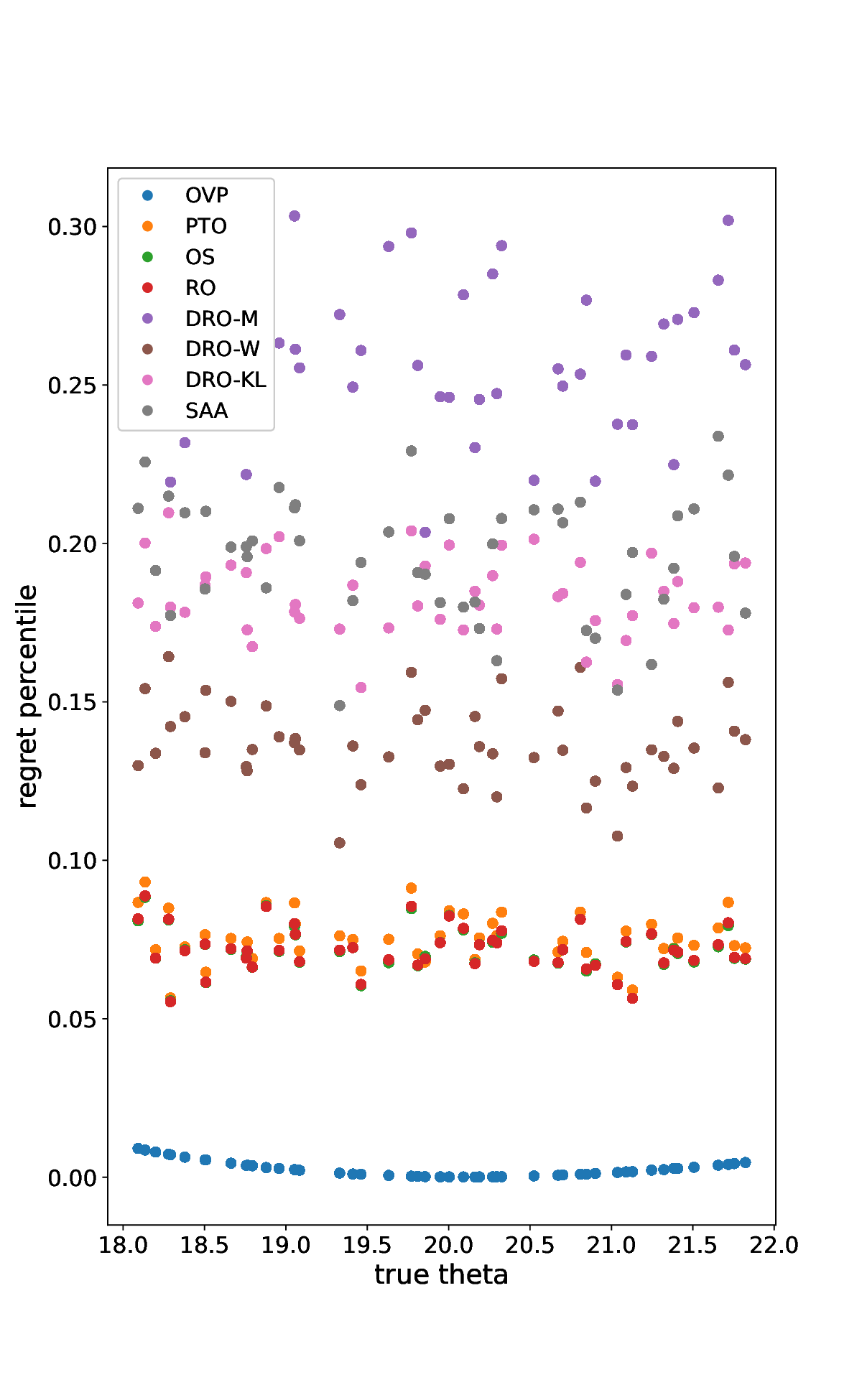} }}
    \vspace{-0.8cm}
    \caption{Average profit (top) and ex-ante percentage regret (bottom) of different models over the range of $\theta$ in the localizations $u\sim\mathcal{N}(20,2)$ (left) and  $u\sim\mathcal{U}[18,22]$ (right) }\label{fig.localization}%
    \vspace{-0.3cm}
\end{figure}

\begin{figure}[h] 
    \centering
    \vspace{-1cm}
   {{\includegraphics[height=10cm]{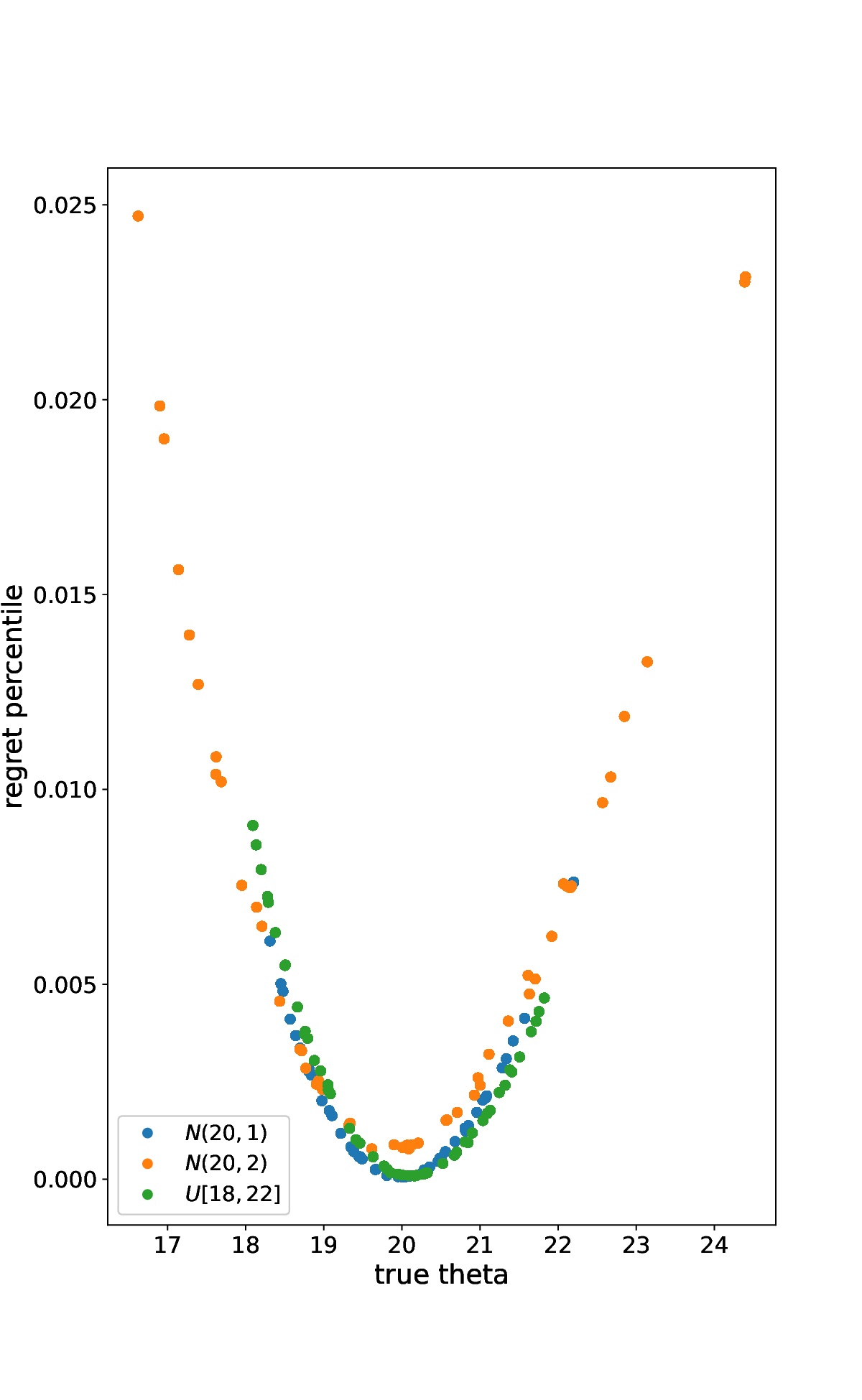}}}
   {{\includegraphics[height=10cm]{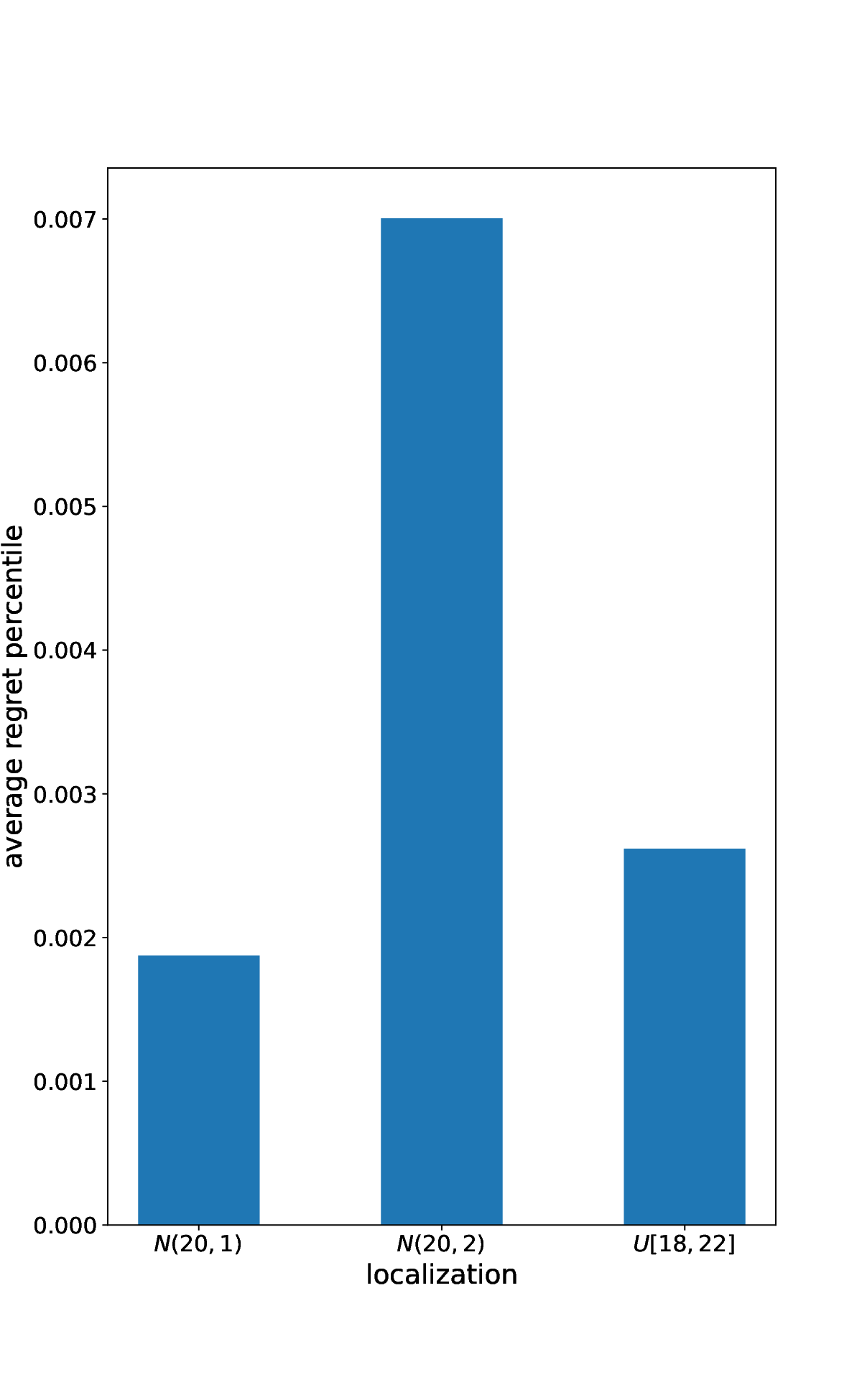} }}
    \vspace{-0.7cm}
    \caption{Average profit (left) and ex-ante percentage regret (right) of OVP for three localizations }\label{fig.ovp.localization}%
    \vspace{-0.5cm}
\end{figure}

In Figure \ref{fig.localization}, we consider localizations $u\sim\mathcal{N}(20, 2)$ (left) and $u\sim U[18, 22]$ (right). The same trends hold. In Figure \ref{fig.ovp.localization}, we compare the OVP solutions obtained over the three different localizations. OVP solutions for different localizations are pareto optimal, \eg, we are unable to tell if the solution for localization $u\sim\mathcal{N}(20,1)$ outperforms $u\sim\mathcal{U}[18,22]$. However, the wider the variance of $u$, \ie, the range of $\theta$ OVP accommodates, the poorer it performs on each specific $\theta$, \eg, $u\sim\mathcal{N}(20,2)$ leads to a worse regret than $u\sim\mathcal{N}(20,1)$. This is the specificity-sensitivity trade-off. 

\subsection{Experiment 2: With misspecification}\label{subsec.misspec}
\begin{figure}[h] 
    \centering
    \vspace{-0.5cm}
    {{\includegraphics[height=10cm]{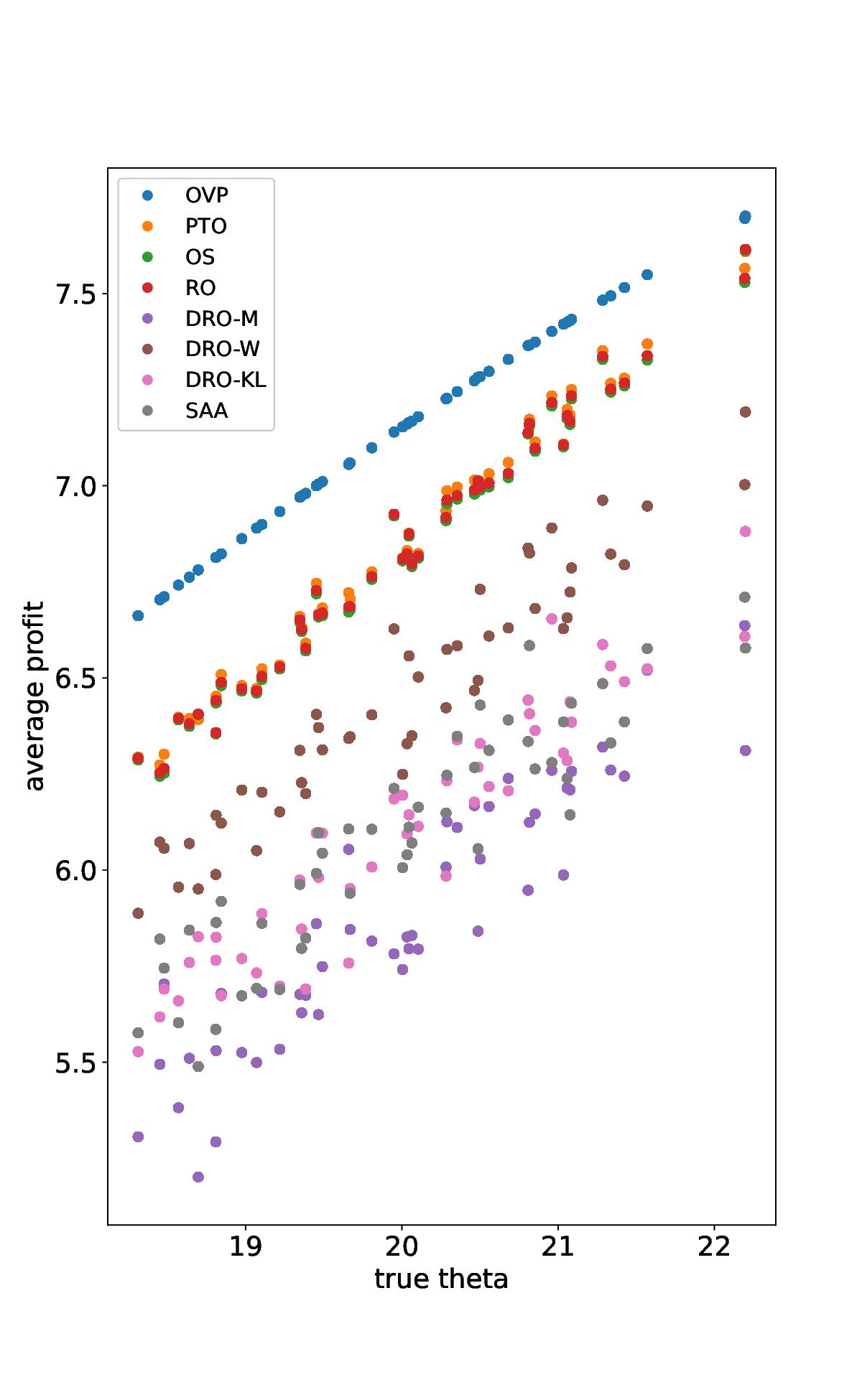}}}
   {{\includegraphics[height=10cm]{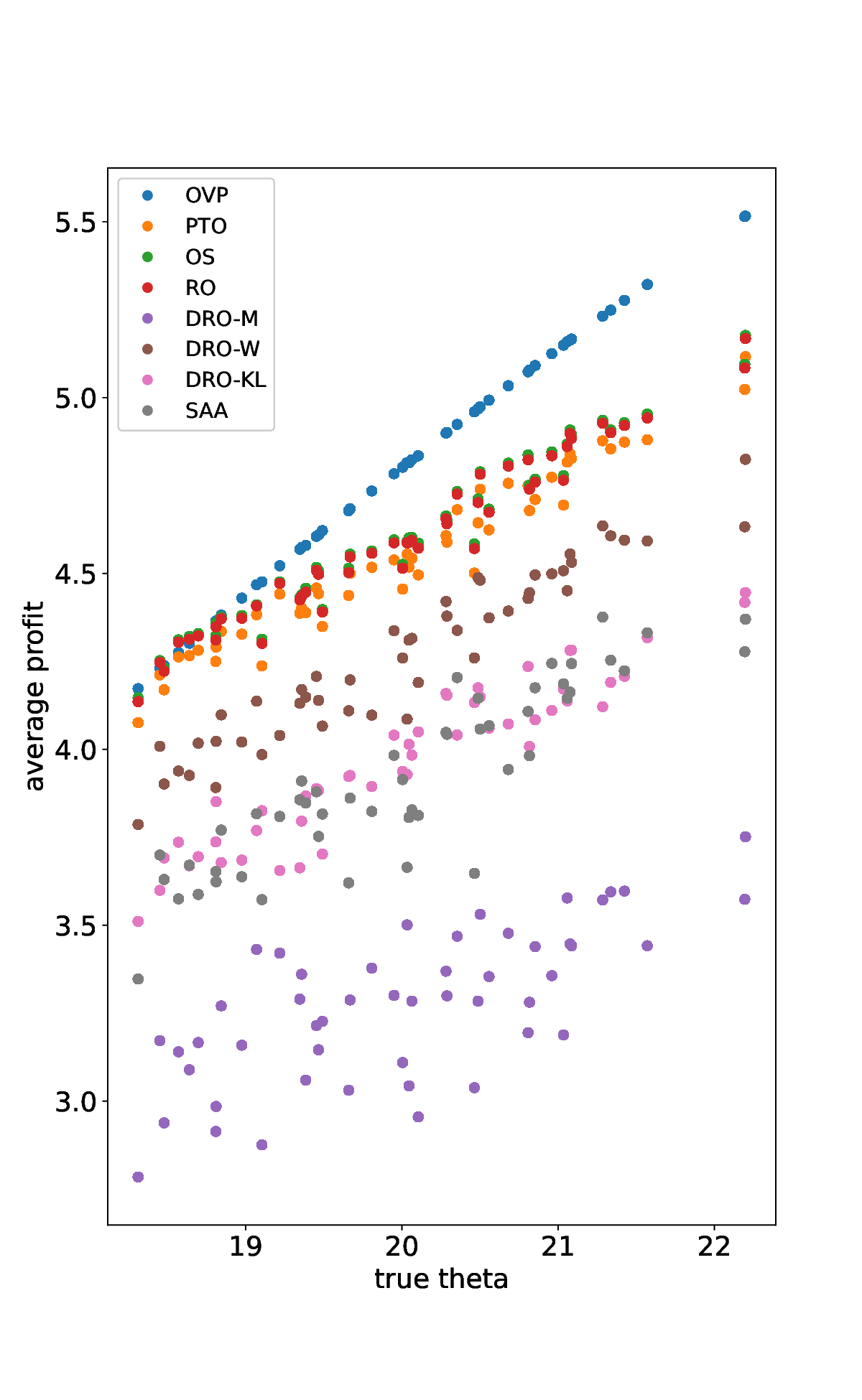} }}%
    \vspace{-0.8cm}
    \caption{Average profit for true demand distribution $d\sim\mbox{Gamma}(1.15, \theta)$ (left) and  $d\sim\mbox{Gamma}(0.85, \theta)$ (right)}\label{fig.profit_mis}%
    \vspace{-0.3cm}
\end{figure}

Consider the case where the true distribution was a Gamma distribution, but the assumed family is the exponential distribution. This implies misspecification. Figure \ref{fig.profit_mis} shows the average profit when demand is distributed by $d\sim\mbox{Gamma}(1.15, \theta)$ and $d\sim\mbox{Gamma}(0.85, \theta)$ respectively. Note that $\mbox{Gamma}(1, \theta) \sim \mbox{Exp}(\theta)$. Models assuming exponential demand are affected (PTO, OS, RO and OVP), whereas data-driven models (SAA and DRO) are more immune, though they have yet to fully close the gap. To fix misspecification, one could solve OVP with a larger hypothesis family, \eg, the Gamma distribution, with localization centred around $1$ for the shape parameter. However, one needs to pay the price of the specificity-sensitivity trade-off.

\section{Conclusions}\label{sec.conclude}

We realized a means of out-of-sample optimal solutions for data-driven optimization. As well, this opens a new chapter on data-driven optimization in regards to misspecification and the selection of the localization. Our work also opens the door to the tantalizing possibilities of out-of-sample optimal end-to-end-learning and bayesian optimization.

\medskip
\bibliographystyle{ormsv080}
\bibliography{bib}

\newpage

\begin{appendices}

\section{Further details of the numerical results}\label{append.num}

In this Appendix, we provide further details on all numerical experiments performed.

\subsection{Computational strategies}\label{append.codes}

In most cases, since $\phi$ is an integral, it is usually differentiable, even if $C$ is not (\emph{e.g.}, the newsvendor problem). In this case, one might draw a sample $\mathcal{U}^{M}:= \{\theta_m\}_{m=1}^{M}$ under the localization $u$, which is known, and directly solve the bisection search problem on
\begin{equation*}
    \sum\limits_{m=1}^{M} \frac{\partial \phi}{\partial q}(q, \theta_{m}) g_1(\hat\theta, \theta_m) = 0.
\end{equation*}

In Algorithm \ref{algo.news}, we present a sample bisection search algorithm for solving for OVP given a particular data set (in the algorithm, only $\hat\theta$ is required) and a localization, approximated by a sample $\mathcal{U}$, for the newsvendor problem with exponential demand.

\begin{algorithm}[!ht]
	\caption{Solving OVP for newsvendor problem with exponential demand} \label{algo.news}
	\begin{algorithmic}[1]
        \Require $p, c, N, \mathcal{U}, \epsilon$. Let $a < b$ given with $b$ sufficiently large. Let $S$ sufficiently large, \emph{e.g.}, $S$ is the mean of $\mathcal{U}$ raised to the power of $N$.
        \Function{Search\_obj}{$q, \hat\theta$}:         
        \State $\text{Search\_obj}  \displaystyle \gets S \cdot \sum\limits_{\theta_m \in \mathcal{U}} \frac{c - pe^{-q/\theta_m}}{\theta_m^N}e^{-N \hat\theta / \theta_m}$
        \EndFunction
        \Function{OVPsolve}{$\hat\theta$}:
        \State $u \gets b$
        \State $l \gets a$
        \While{$u-l> \epsilon$}
        \State $q \gets (u+l)/2$ 
        \If{$\text{Search\_obj}(q, \hat\theta) = 0$}
          \Return $q$
        \ElsIf{$\text{Search\_obj}(q, \hat\theta) \times \text{Search\_obj}(l, \hat\theta)< 0$} 
        \State $b \gets q$ 
        \Else
        \State $a \gets q$
        \EndIf    
        \EndWhile  \\
        \Return $q$
        \EndFunction
	\end{algorithmic} 
\end{algorithm}

Note that it is likely that $g_1$ is a very small number, because it was originally the density of a high-dimensional integral, but is now defined only on a subspace subtended by $\hat\theta$. In the case of the newsvendor model, it would involve large divisions by $\theta_m^N$ and $e^{-N \hat\theta / \theta_m}$. In view of this, the large constant $S$ is incorporated in the Algorithm to circumvent numerical stability issues.

In the event that $C$ is convex, but not differentiable, one can approximate $\beth$ with the family of datasets $\mathcal{F}^{\bar{N}}(\theta_m)=\{d_{m,l}\}_{l=1}^{\bar{N}}$ under the distribution of $f(\cdot, \theta_m)$ and implement a convex optimization on an approximate $\beth$. Once again, as both $u$ and $f$ are known, one may draw as many samples as one desires to approximate $\beth$ to arbitrary accuracy (in exchange for computational efficiency):
\begin{equation*}
    \min_q \quad \sum\limits_{m=1}^{M} \sum\limits_{l=1}^{\bar{N}}  C(q, d_{m,l}) g_1(\hat\theta, \theta_m).
\end{equation*}

In Algorithm \ref{algo.news_compare}, we present a sample pseudocode for evaluating the out-of-sample performance for a prescriptive solution, such as OVP, for the newsvendor problem with exponential demand.

\begin{algorithm}[!ht]
	\caption{Evaluating performance of OVP and benchmarks for exponential newsvendor} \label{algo.news_compare}
	\begin{algorithmic}[1]
        \Require $p, c, N, \mathcal{U}, \bar{N}$. Let $M = |\mathcal{U}|$. Functions from Algorithm \ref{algo.news}. Let $\mbox{Solve}$ be some function for obtaining the policy, \emph{e.g.}, `OVPsolve'.
        \Function{True\_cost}{$q, \theta$}:
        \State $\text{True\_cost}  \gets qc - p\theta\left( 1- e^{-q\theta} \right)$
        \EndFunction
        \For {$\theta_m$ in $\mathcal{U}$}
        \For {$i=1$ to $\bar{N}$}
        \State Sample $d_1, \ldots, d_N$ from $f(d,\theta_m)$.
        \State $\displaystyle \hat\theta \gets \frac{1}{N} \sum_{n=1}^N d_n$
        \State $q \gets \text{Solve}(\hat\theta)$ 
        \State $\phi[i] \gets \text{True\_cost}(q, \theta_m)$ 
        \EndFor
        \State $\displaystyle \Phi[\theta_m] \gets \frac{1}{\bar{N}} \sum_{i=1}^{\bar{N}} \phi[i]$
        \EndFor \\
        \Return $ \displaystyle \frac{1}{M}\sum_{m=1}^{M} \Phi[\theta_m] $
	\end{algorithmic} 
\end{algorithm}

\subsection{Mathematical formulations for benchmark models}\label{append.benchmarks}

Here we clearly state the formulations for the benchmark models used in our numerical experiments, specific to the case of the newsvendor model. 

\smallskip
\noindent\ul{Operational statistic solution}: Here, the solution $q_{\scriptscriptstyle{\text{OS}}} = \alpha \hat\theta(\bm y)$ is considered, where $\alpha$ is to be determined optimally via:
\begin{equation*}
    \min\limits_{\alpha} \mathbb{E}_{\bm y}[ \alpha \hat\theta(\bm y) - p\theta( 1- e^{-\alpha \hat\theta(\bm y) / \theta}) \,|\, \theta ].
\end{equation*}
This solution turns out to be independent of $\theta$ and has the closed form  $\alpha = N\left((\frac{p}{c})^{1/(N+1)}-1\right)$, where $N$ is the number of data points

\begin{equation*}
q^* = \theta \log(p/c)
\end{equation*}

\smallskip
\noindent\ul{Robust optimization model with $\theta$ as the uncertainty}: The robust newsvendor problem can be expressed as 
\begin{equation*}
     \min_{q\geq 0} \sup_{\theta \in \Omega }\left\{   \phi(q, \theta) \right\}   = \min_{q\geq 0} \sup_{\theta \in \bar\Omega }\left\{  qc - p\theta\left(1- e^{-q/\theta}\right)  \right\},
\end{equation*}
for some uncertainty set $\Omega$. In the context of our setting, $\theta$ is the mean-parameter of the exponential distribution and is one-dimensional. Requiring $\Omega$ to be closed and convex, $\Omega$ will essentially be a closed interval -- $\Omega = [\underbar{$\omega$}, \bar\omega]$ containing $\hat\theta$. We specify these bounds as a fraction of $\hat\theta$, specifically, $\Omega = [0.95\,\hat\theta, 1.05\,\hat\theta]$.

Note that $ \phi(q, \theta)$ is jointly convex in $q$ and $\theta$. Hence the worst-case $\theta$ must belong to the boundary of $\Theta$, and the problem simplifies to 
\begin{equation*}
     \min_{q\geq 0} \max\left\{  qc - p\underbar{$\omega$}\left(1- e^{-q/\underbar{$\omega$}}\right) , qc - p\bar\omega\left(1- e^{-q/\bar\omega}\right)  \right\} = \min_{q\geq 0} \left\{  qc - p\underbar{$\omega$}\left(1- e^{-q/\underbar{$\omega$}}\right)   \right\}.
\end{equation*}
The equation above is due to the fact that  $\frac{\partial}{\partial \theta} (qc - p\theta(1-e^{-q/\theta)}) \leq 0$ for all $q\geq 0$ and $\theta >0$. Thus, this problem has closed form solution $q^* = \underbar{$\omega$} \log (p/c)$. In other words, `the worst-case scenario is always when the demand is smaller than expected'.

\smallskip
\noindent\ul{Distributionally robust optimization model with moment uncertainty}: The distributionally robust optimization formulation for the newsvendor problem can be written as: 
    \begin{equation}\label{eq.dro}
        \min_{q\geq 0} \sup_{\mathbb{P}\in \mathcal{P}}\left\{  \mathbb{E}_{\mathbb{P}}\left[ qc -  p \min\{\tilde{d}, q\}  \right]  \right\},
    \end{equation}
  for some ambiguity set $\mathcal{P}$. In the case of moments uncertainty, $\mathcal{P}$ is 
    \begin{equation*}
         \mathcal{P}_m = \left\{ \mathbb{P}\in \mathcal{P}(\mathbb{R}_+)~\left|
         \begin{array}{ll}
              & \tilde{d}\sim \mathbb{P} \\
              & \mathbb{E}_{\mathbb{P}} \left[  \tilde{d} \right] = \hat{d} \\
               & \mathbb{E}_{\mathbb{P}} \left[  \left (\tilde{d} - \hat{d} \right)^2\right] \leq  \hat{\sigma}^2 
         \end{array}
         \right\} \right. ,
    \end{equation*}
    with robust counterpart
    \begin{equation*}
        \begin{array}{rll}
           \inf  &  \lambda \hat{d} + \beta \hat{\sigma}^2 + \gamma \\
           \mbox{s.t.}   &  \lambda d + \beta (d-\hat{d})^2 + \gamma \geq qc - p\min\{d,q\} ~~ \forall d\geq 0, \\
                         & \beta \geq 0,\, \lambda,\gamma ~\text{free}.
        \end{array}
    \end{equation*}

\smallskip
\noindent\ul{Distributionally robust optimization model with Wasserstein uncertainty set}: 
We instead consider the ambiguity set $ \mathcal{P}_W^r $ in \eqref{eq.dro},
    \begin{equation*}
         \mathcal{P}_W^r = \left\{ \mathbb{P}\in \mathcal{P}(\mathbb{R}_+)~\left|
         \begin{array}{ll}
              & \tilde{d}\sim \mathbb{P} \\
              & \Delta_W(\mathbb{P},\hat{\mathbb{P}}) \leq r
         \end{array}
         \right\} \right.,
    \end{equation*}
    where $\Delta_W(\mathbb{P},\hat{\mathbb{P}}) $ is the Wasserstein distance defined on some norm $|| \cdot ||$, and the empirical distribution $\hat{\mathbb{P}}$ is given by $\hat{\mathbb{P}}[\tilde{d} = \hat{d}_i] = 1/N$, for all $i=1,\ldots,N$. Its robust counterpart is
    \begin{equation*}
        \begin{array}{rll}
            \inf & \displaystyle \lambda r + \frac{1}{N}\sum_{i=1}^N \beta_i \\
            \mbox{s.t.} & \lambda \| d - \hat{d}_i\| + \beta_i \geq qc -p\min\{d,q\} ~~ \forall d\geq 0,\, i=1,\ldots,N,\\
            & \lambda \geq 0,\, \beta_i~\text{free}, i=1,\ldots,N.
        \end{array}
    \end{equation*}
In our case, we shall just consider the $L_1$-norm, and this problem can be easily solved as a linear program. We searched for the radius $r$ via cross-validation and how this is done is explained in Appendix \ref{append.cross}.

\smallskip
\noindent\ul{Distributionally robust optimization model with Kullback-Leibler uncertainty set}: We instead consider the ambiguity set $ \mathcal{P}_{KL}^r $ in \eqref{eq.dro},
    \begin{equation*}
         \mathcal{P}_{KL}^r = \left\{ \mathbb{P}\in \mathcal{P}(\mathbb{R}_+)~\left|
         \begin{array}{ll}
              & \tilde{d}\sim \mathbb{P} \\
              & \Delta_{KL}(\mathbb{P},\hat{\mathbb{P}}) \leq r
         \end{array}
         \right\} \right.,
    \end{equation*}
    and $\Delta_{KL}(\mathbb{P},\hat{\mathbb{P}}) $ is the Kullback-Leibler divergence. Let $\hat{p}_i = \hat{\mathbb{P}}[\tilde d = \hat{d}_i]$, for all $i=1,\ldots,N$.  The corresponding robust counterpart is 
\begin{equation*}
    \begin{array}{rll}
        \inf & \displaystyle \lambda r + \sum_{i=1}^N \beta_i \hat{p}_i + \gamma  \\
        \mbox{s.t.} &  \lambda \log\left( \lambda/\beta_i\right) + qc - p\min\{\hat{d}_i,q\} \leq \lambda +\gamma~~\forall i = 1, \ldots,N,\\
                    & \lambda \geq 0,\, \gamma~\text{free},\, \beta_i\geq 0,~i = 1,\ldots, N.
    \end{array}
\end{equation*}
We follow the convention that $0\log(0/t)=0$ if $t\geq 0$. This problem can be expressed as an exponential cone program, and we solved it using the exponential cone solver in MOSEK. We searched for the radius $r$ via cross-validation and how this is done is explained in Appendix \ref{append.cross}.

\smallskip
\noindent\ul{DRO model with moment uncertainty on assumed Exponential demand}: The distributionally robust optimization model assuming uncertainty in the unknown parameter $\theta$ can be expressed as 
\begin{equation*}
     \min_{q\geq 0} \sup_{\mathbb{P} \in \mathcal{P}_m^e }\left\{ \mathbb{E}_{\mathbb{P}} [\phi(q, \theta)] \right\}   = \min_{q\geq 0} \sup_{\mathbb{P} \in \mathcal{P}_m^e }\left\{  \mathbb{E}_{\mathbb{P}}\left[ qc - p\tilde{\theta}\left(1- e^{-q/\tilde{\theta}} \right) \right]  \right\},
\end{equation*}
where
\begin{equation*}
         \mathcal{P}_m^e = \left\{ \mathbb{P}\in \mathcal{P}(\mathbb{R}_+)~\left|
         \begin{array}{ll}
              & \tilde{\theta}\sim \mathbb{P} \\
              & \mathbb{E}_{\mathbb{P}} \left[  \tilde{\theta} \right] = \hat{\theta} \\
               & \mathbb{E}_{\mathbb{P}} \left[  \left (\tilde{\theta} - \hat{\theta} \right)^2\right] \leq  s^2 
         \end{array}
         \right\} \right. ,
    \end{equation*}
    and $s^2$ is the unbiased estimator for the variance of the sample mean $\hat\theta$. Unfortunately, this model is non-convex. As such, we do not consider it in our simulations.

\subsection{Simulation set-up}\label{append.setup}

The simulation set-up is standardized for all of the numerical experiments conducted in Section \ref{sec.illust}. The same settings are used in the illustration in Example \ref{example.quad}, except $N=5$ is used instead.

\smallskip
\noindent\ul{Parameters}: The selling price is set at $p=2$ and the cost price is set at $c=1$. Sensitivity analysis on the prices were not conducted, but we had done a quick check on a different set of prices to realize that the key findings and insights had not changed. For the large constant $S$ in Algorithm \ref{algo.news_compare}, we had used $\check\theta^N e^{N\hat\theta / \grave{\theta}}$, where $\check\theta = \min \theta_m$ and $\grave{\theta} = \max \theta_m$. The upper and lower bounds for the bisection search in Algorithm \ref{algo.news} are given as $a = 0$ and $b = 2\hat\theta$, with automatic subsequent relaxation of $b$ if it is initially tested to be of the same sign as the solution generated slightly above $a$.

\smallskip
\noindent\ul{Generating datasets}: We chose $M=50$, that is, the number of samples to draw from the localization, $N = 10$, that is, the number of observed demand data points per data set, and $\bar{N} = 200$, that is the number of different data sets we re-sampled in order to reflect the sampling distribution. We chose $N = 10$ as differences between the models are sufficiently pronounced for clear comparisons. Findings are consistent even if $N$ is increased. When $N \geq 20$ roughly, the number of data points is sufficient to estimate $\hat\theta$ to relatively high accuracy, thus all solutions begin to converge to the OVP solution. Notice that none of the benchmark models depend on $M$ and $\bar{N}$. We chose $\bar{N}$ to be sufficiently large as we realized that there is reasonable amount of variation in performance across data sets (which also arises because the variation in $\hat\theta$ is large when $N$ is small). In particular, we chose $\bar{N}$ to be large enough so that OS performs better than PTO on almost every $\theta$, as is theoretically known. This would be a good indication that the variations have been sufficiently averaged away, and that happens roughly around $\bar{N} \geq 200$, which is what we have chosen. When $M=20$, OVP already exhibits very clear distinction against the other benchmark models, but we chose $M=50$ so that we would obtain a better spread of test true $\theta$'s, and to give ample chance for outliers to occur to test the generalizability of the models. Notice that for the localization, we had used different sets of ($M=50$) points for computing the OVP solution in Algorithm \ref{algo.news} versus the out-of-sample performance in Algorithm \ref{algo.news_compare}, so as to ensure that OVP would be able to generalize regardless of the sample used for the localization. 

Within each experiment, for all the benchmarks, we used the same set of samples for the localization, and data sets. Only the data sets used for the cross-validation for the Wasserstein and KL divergence radii differ (discussed later in Appendix \ref{append.cross}). 

\subsection{Cross-validation for Wasserstein and KL radii}\label{append.cross}

In our simulations, we aim to obtain the best parameter for the radii for the Wasserstein and KL divergence uncertainty sets to illustrate the limits of these approaches. In other words, the results presented for these two models are already conditional on having chosen the best possible radii $r$. As such, the procedure here would be possibly considered `counterfactual' since in each instance, the decision-maker would only possess the training data set, and can only conduct cross-validation on that training data set. This is not even reasonable in some cases, such as in this simulation experiment where the number of data points in each data set is $N=10$. 

Nonetheless, notice that because generalizability is part of what we are interested in, specifically that the model works well over a range of true $\theta$'s, we need to apply the same radius for different $\theta$. We randomly generate $20$ samples of $\theta$'s from the localization $u$. Based on each $\theta$, we generate a data set of size $N=10$ and compute the average performance of these data-driven models over these data sets. The grid search for the optimal radius $r$ is conducted over the range of $[0.0001, 5]$. To evaluate the performance of each solution, both in the no misspecification and misspecification experiments, we calculate the percentage gap against the true ex-ante oracle. Once again, this is counterfactual, but for the purposes of comparison, we have done so for to maximize the potential of the benchmarks. As discussed in \S\ref{sec.illust}, the OVP solution incurs almost zero regret for the no misspecification case, and represents the limit of performance attainable with learning (\ie, without perfect information). We then select the best radius $r$ that gives the smallest average percentage gap. 

While it is possible to compute the cross-validation performance for every $\theta$ and every data set for that given $\theta$, this would take a lot of time, and instead, this sampling procedure is adopted.

The cross-validated radius is separately computed for each localization. Figures \ref{fig.cv_norm_1} to \ref{fig.cv_mis_85} below show the average percentage gaps for the different radius over the search range. 

\begin{figure}[!h]
    \centering
    \vspace{-0.5cm}
    \includegraphics[width = 10cm]{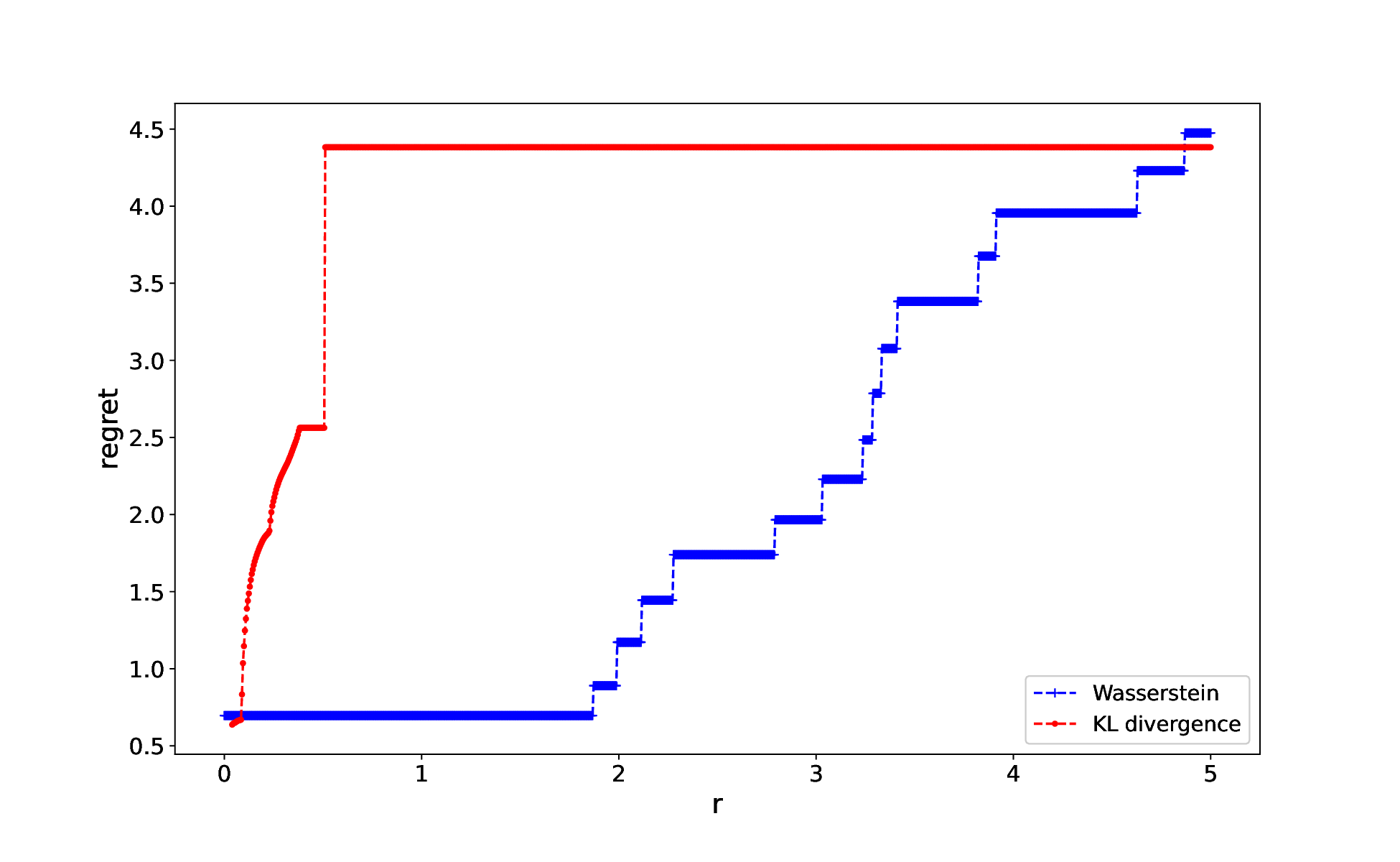}
    \caption{Cross validation for localization $\mathcal{N}(20,1)$}
    \label{fig.cv_norm_1}
\end{figure} 

\begin{figure}[!h]
    \centering
    \vspace{-0.5cm}
    \includegraphics[width = 10cm]{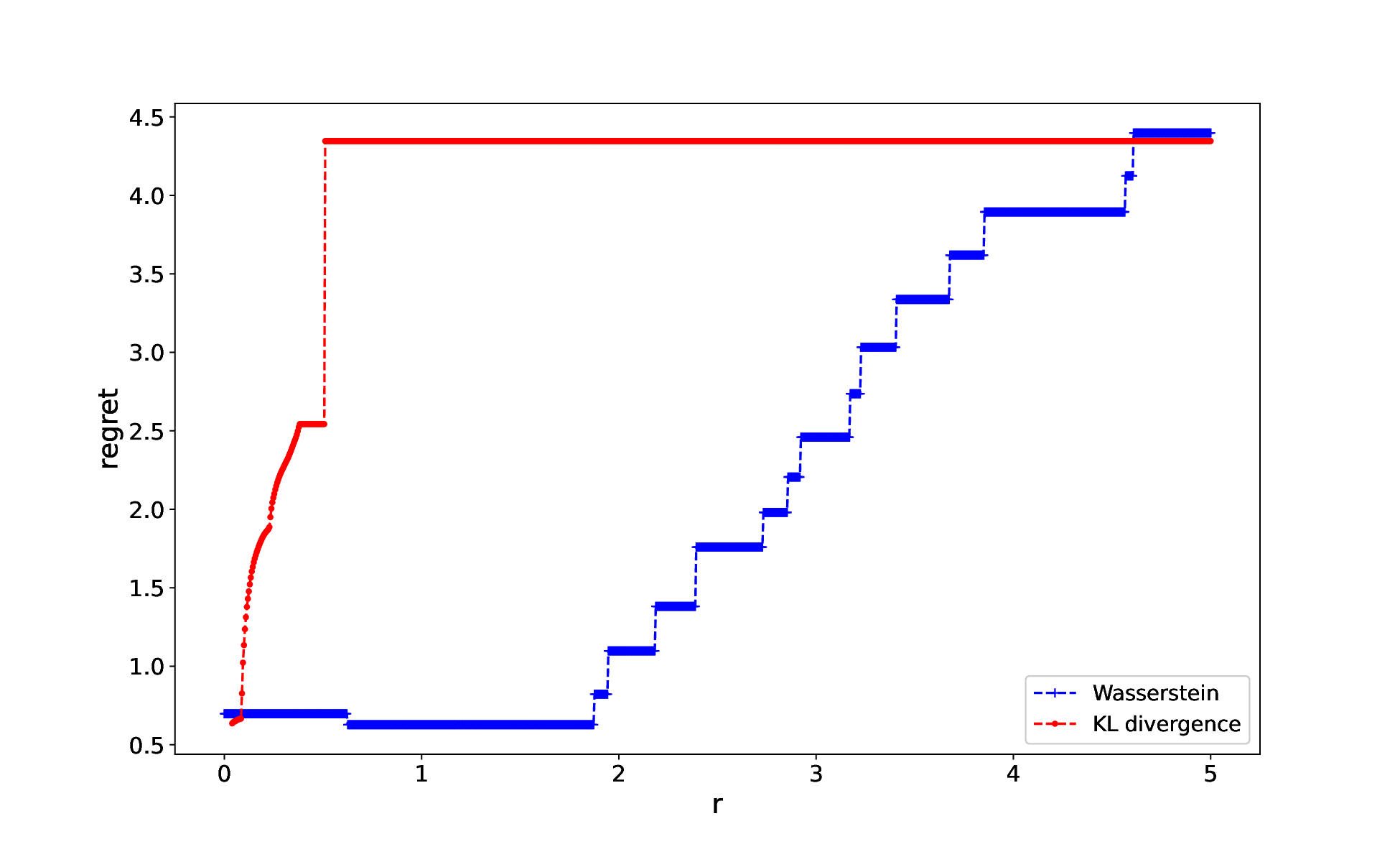}
    \caption{Cross validation for localization $\mathcal{N}(20,2)$}
    \label{fig.cv_norm_2}
\end{figure} 

\begin{figure}[!h]
    \centering
    \vspace{-0.5cm}
    \includegraphics[width = 10cm]{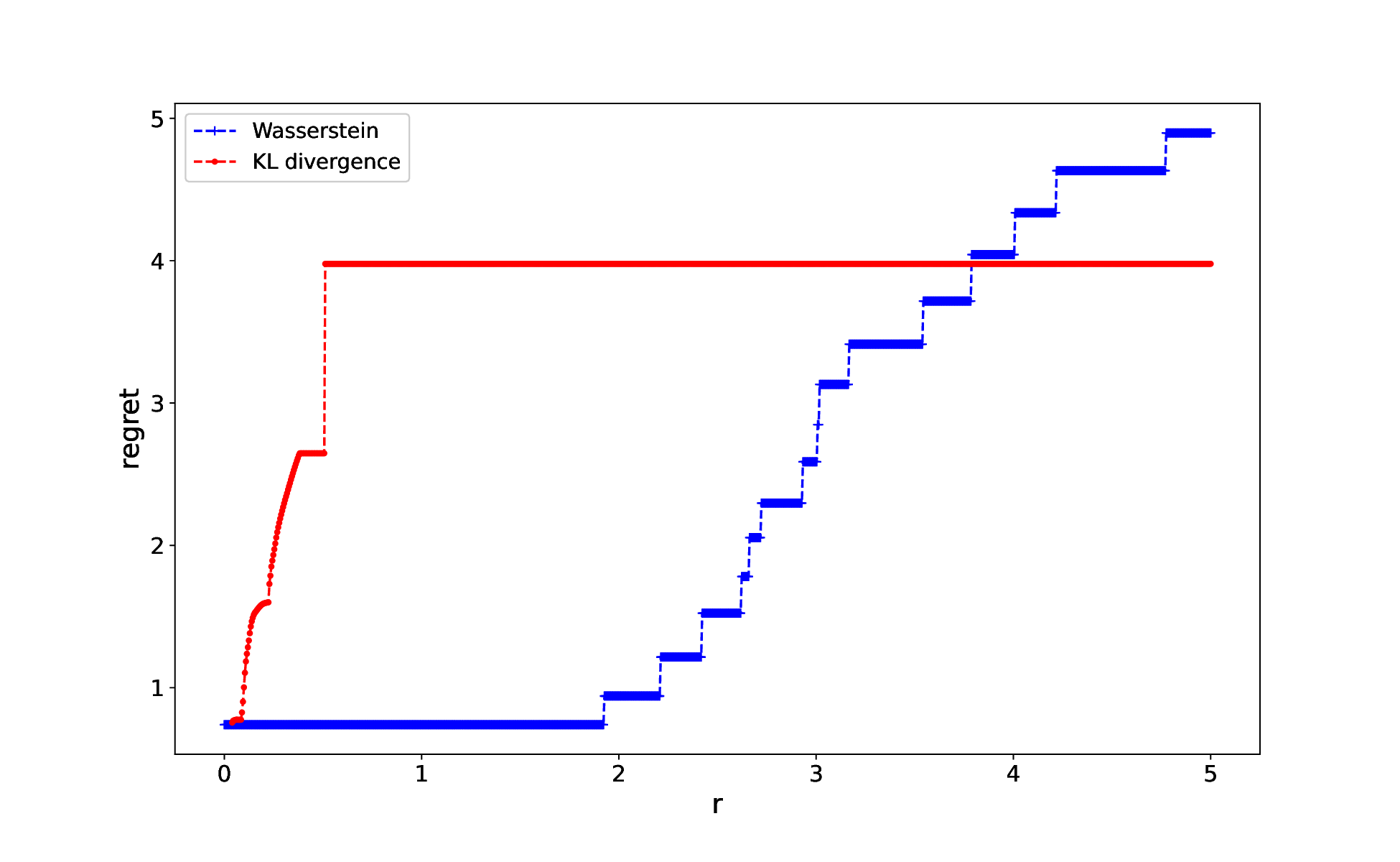}
    \caption{Cross validation for localization $\mathcal{U}(18,22)$}
    \label{fig.cv_uniform}
\end{figure} 

\begin{figure}[!h]
    \centering
    \vspace{-0.5cm}
    \includegraphics[width = 10cm]{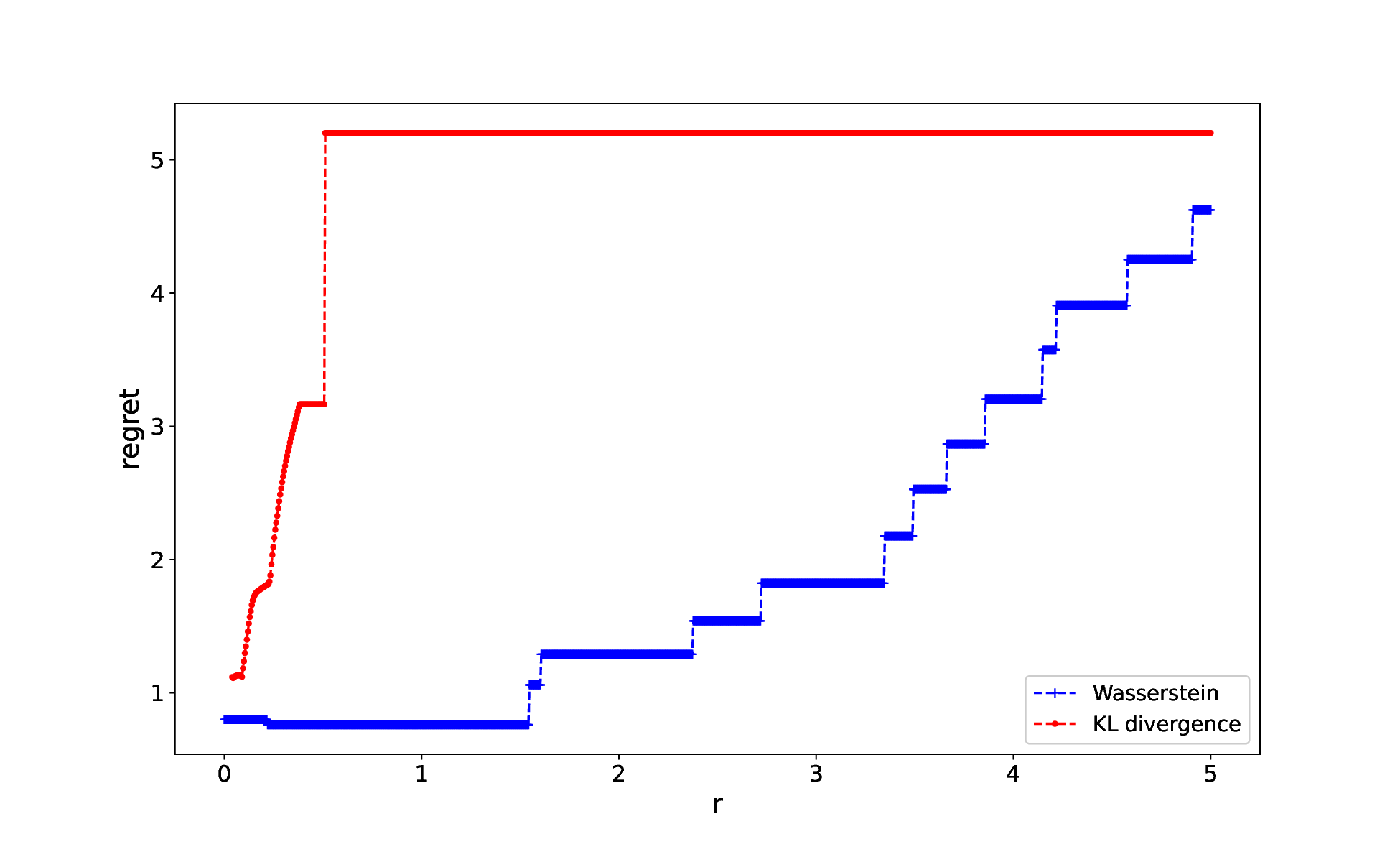}
    \caption{Cross validation for misspecificiation $d\sim\mbox{Gamma}(1.15, \theta)$ }
    \label{fig.cv_mis_15}
\end{figure}

\begin{figure}[!h]
    \centering
    \vspace{-0.5cm}
    \includegraphics[width = 10cm]{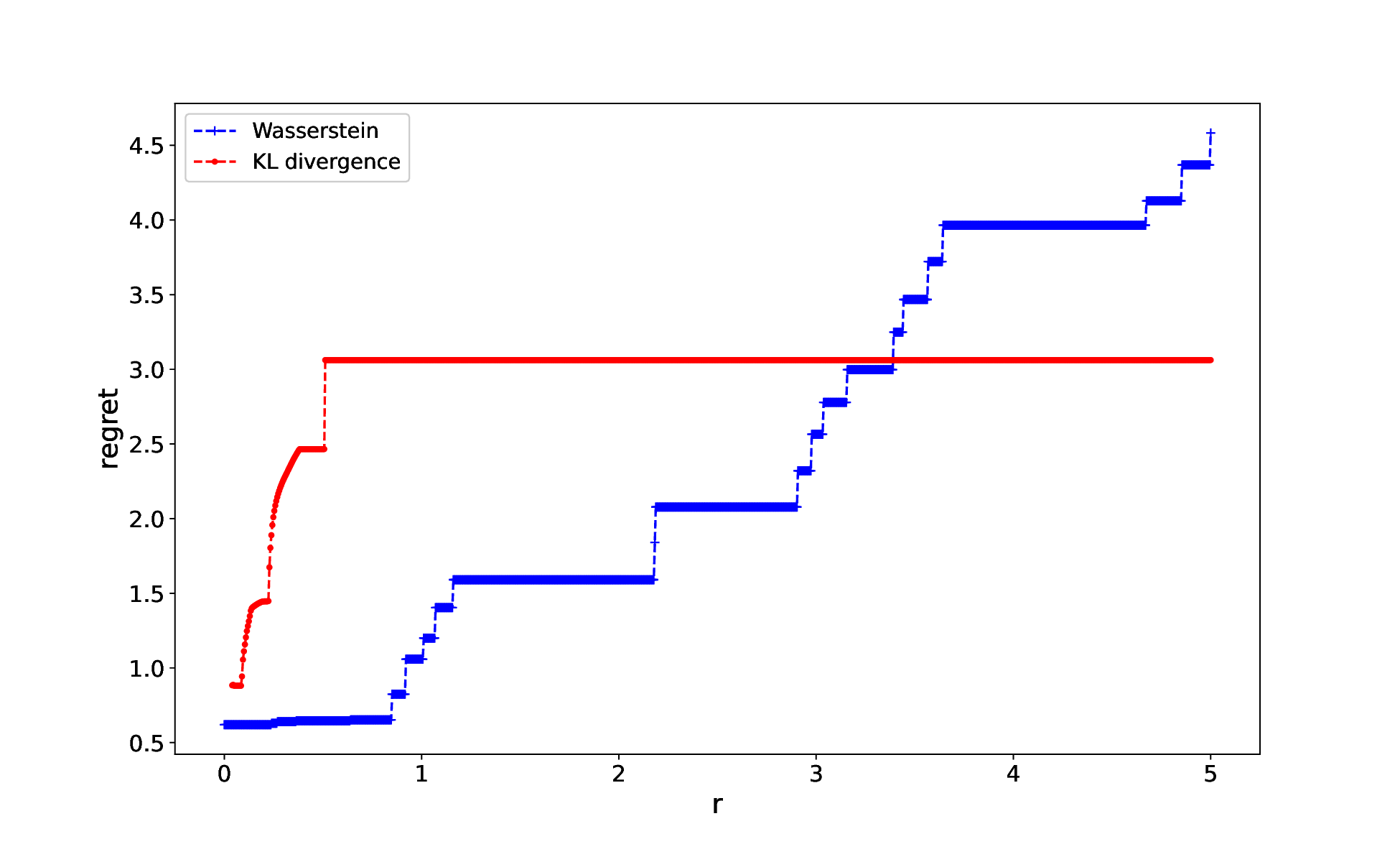}
    \caption{Cross validation for misspecification $d\sim\mbox{Gamma}(0.85, \theta)$ }
    \label{fig.cv_mis_85}
\end{figure}

\end{appendices}

\end{document}